\documentclass[12pt,draftcls,onecolumn]{IEEEtran}
\usepackage{dsfont}
\usepackage{amssymb}
\usepackage{amsmath}
\usepackage[dvips]{graphicx}
\usepackage{subfigure}

\newtheorem{thm}{Theorem}
\newtheorem{lem}{Lemma}
\newtheorem{rmk}{Remark}
\newtheorem{defn}{Definition}
\newtheorem{prop}{Proposition}
\def\<{\langle}
\def\>{\rangle}

\def \w {\omega}
\def \t {\theta}
\def \a {\alpha}
\def \b {\beta}

\def\bR {\mathbb{R}}

\begin{document}
\title{Control and Synchronization of Neuron Ensembles}
\author{Jr-Shin~Li,~%\IEEEmembership{Member,~IEEE}
        Isuru Dasanayake,~%\IEEEmembership{Student Member,~IEEE,}
        and~Justin~Ruths% <-this % stops a space
%\thanks{%Manuscript received January 20, 2002; revised November 18, 2002.
%		This work was supported in part by the National Science Foundation under the Career Award 0747877 and the Air Force Office of
%		Scientific Research Young Investigator Award FA9550-10-1-0146.}
        % This work was supported by the NSF SGER, \#0744090, and CAREER, \#0747877, grants.}% <-this % stops a space
\thanks{J.-S. Li is with the Department of Electrical and Systems Engineering, Washington University, St. Louis, MO 63130 USA (e-mail: jsli@seas.wustl.edu).}
\thanks{I. Dasanayake is with the Department of Electrical and Systems Engineering, Washington University, St. Louis, MO 63130 USA (e-mail: dasanayakei@seas.wustl.edu).}
\thanks{J. Ruths is with the Engineering Systems \& Design Pillar, Singapore University of Technology \& Design, Singapore (e-mail: justinruths@sutd.edu.sg).}}

\maketitle

% ================================== Abstract =================================
\begin{abstract}
	Synchronization of oscillations is a phenomenon prevalent in natural, social, and engineering systems. Controlling synchronization of oscillating systems is motivated by a wide range of applications from neurological treatment of Parkinson's disease to the design of neurocomputers. In this article, we study the control of an ensemble of uncoupled neuron oscillators described by phase models. We examine controllability of such a neuron ensemble for various phase models and, furthermore, study the related optimal control problems. In particular, by employing Pontryagin's maximum principle, we analytically derive optimal controls for spiking single- and two-neuron systems, and analyze the applicability of the latter to an ensemble system. Finally, we present a robust computational method for optimal control of spiking neurons based on pseudospectral approximations. The methodology developed here is universal to the control of general nonlinear phase oscillators.
\end{abstract}


\IEEEpeerreviewmaketitle

%%%%%%%%%%%%%%%%%%%%%%%%%%%%%%%%%%%%%%%%%%%%%%%%%%%%%%%%%%%%%%%%%%%%%%%%
\section{Introduction}
\label{intro}
Natural and engineered systems that consist of ensembles of isolated or interacting nonlinear dynamical components have reached levels of complexity that are beyond human comprehension. These complex systems often require an optimal hierarchical organization and dynamical structure, such as synchrony, for normal operation. The synchronization of oscillating systems is an important and extensively studied phenomenon in science and engineering \cite{Strogatz01}. Examples include neural circuitry in the brain \cite{Uhlhaas06}, sleep cycles and metabolic chemical reaction systems in biology \cite{Hanson78,Mirollo90,Ermentrout84,Nishikawa08}, semiconductor lasers in physics \cite{Fischer00}, and vibrating systems in mechanical engineering \cite{Blekhman88}. Such systems, moreover, are often tremendously large in scale, which poses serious theoretical and computational challenges to model, guide, control, or optimize them. Developing optimal external waveforms or forcing signals that steer complex systems to desired dynamical conditions is of fundamental and practical importance \cite{Harada10,Kiss07}. For example, in neuroscience devising minimum-power external stimuli that synchronize or desynchronize a network of coupled or uncoupled neurons is imperative for wide-ranging applications from neurological treatment of Parkinson's disease and epilepsy \cite{Ashwin92, Benabid91, Schiff94} to design of neurocomputers \cite{Hoppensteadt01,Hoppensteadt00}; in biology and chemistry application of optimal waveforms for the entrainment of weakly forced oscillators that maximize the locking range or alternatively minimize power for a given frequency entrainment range \cite{Harada10, Li_DSCC11} is paramount to the time-scale adjustment of the circadian system to light \cite{Winfree80} and of the cardiac system to a pacemaker \cite{Glass91}. 

Mathematical tools are required for describing the complex dynamics of oscillating systems in a manner that is both tractable and flexible in design. A promising approach to constructing simplified yet accurate models that capture essential overall system properties is through the use of phase model reduction, in which an oscillating system with a stable periodic orbit is modeled by an equation in a single variable that represents the phase of oscillation \cite{Winfree80, Kuramoto84}. Phase models have been very effectively used in theoretical, numerical, and more recently experimental studies to analyze the collective behavior of networks of oscillators \cite{Acebron05,Kiss05,Preyer05,Netoff05}. Various phase model-based control theoretic techniques have been proposed to design external inputs that drive oscillators to behave in a desired way or to form certain synchronization patterns. These include multi-linear feedback control methods for controlling individual phase relations between coupled oscillators \cite{Kano10} and phase model-based feedback approaches for efficient control of synchronization patterns in oscillator assemblies \cite{Kiss07,Rusin10,Zhai08}. These synchronization engineering methods, though effective, do not explicitly address optimality in the control design process. More recently, minimum-power periodic controls that entrain an oscillator with an arbitrary phase response curve (PRC) to a desired forcing frequency have been derived using techniques from calculus of variations \cite{Li_DSCC11}. In this work, furthermore, an efficient computational procedure was developed for optimal control synthesis employing Fourier series and Chebyshev polynomials. Minimum-power stimuli with limited amplitude that elicit spikes of a single neuron oscillator at specified times have also been analytically calculated using Pontryagin's maximum principle, where possible neuron spiking range with respect to the bound of the control amplitude has been completely characterized \cite{Li_PRE11, Li_CDC11_Neuron}. In addition, charge-balanced minimum-power controls for spiking a single neuron has been thoroughly studied \cite{Li_IEEE_CB, Nabi09}.

In this paper, we generalize our previous work on optimal control of a single neuron \cite{Li_PRE11, Li_CDC11_Neuron,Li_IEEE_CB} to consider the control and synchronization of a collection of neuron oscillators. In particular, we investigate the fundamental properties and develop optimal controls for the synchronization of such type of large-scale neuron systems. In Section \ref{sec:neuron}, we briefly introduce the phase model for oscillating systems and investigate controllability of an ensemble of uncoupled neurons for various phase models characterized by different baseline dynamics and phase response functions. Then, in Section \ref{sec:oc}, we formulate optimal control of spiking neurons as steering problems and in particular derive minimum-power and time-optimal controls for single- and two-neuron systems. Furthermore, we implement a multidimensional pseudospectral method to find optimal controls for spiking an ensemble of neurons which reinforce and augment our analytic results.

%\textbf{Finally, we apply the derived optimal controls from the phase-reduced models to their original state-space systems to calibrate the relation and explore the approximation limits of the original system by the reduced phase model. Such an important validation is largely lacking in the literature.}

%%%%%%%%%%%%%%%%%%%%%%%%%%%%%%%%%%%%%%%%%%%%%%%%%%%%%%%%%%%%%%%%%%%%%%%%%%%%%%%%%%%%%%
\section{Control of Neuron Oscillators}
\label{sec:neuron}

% ========================================
\subsection{Phase Models}
The dynamics of an oscillator are often described by a set of ordinary differential equations that has a stable periodic orbit. Consider a time-invariant system 
\begin{equation}
    \label{eq:sys1}
    \dot{x}=F(x,u), \quad x(0)=x_0,
\end{equation}
where $ x(t)\in\mathbb{R}^n$ is the state and $u(t)\in\mathbb{R}$ is the control, which has an unforced stable attractive periodic orbit $\gamma(t)=\gamma(t+T)$ homeomorphic to a circle, satisfying $\dot{\gamma}=F(\gamma,0)$, on the periodic orbit $\Gamma=\{y\in\bR^n\,|\,y=\gamma(t),\ \text{for } 0\leq t< T\}\subset\bR^n$. This system of equations can be reduced to a single first order differential equation, which remains valid while the state of the full system stays in a neighborhood of its unforced periodic orbit \cite{Brown04}. This reduction allows us to represent the dynamics of a weakly forced oscillator by a single phase variable that defines the evolution of the oscillation,
\begin{equation}
    \label{eq:phasemodel}
    \frac{d\t}{dt}=f(\theta)+Z(\theta)u(t),
\end{equation}
where $\theta$ is the phase variable, $f$ and $Z$ are real-valued functions, and $u\in\mathcal{U}\subset\mathbb{R}$ is the external stimulus (control) \cite{Brown04, Izhikevich07}. The function $f$ represents the system's baseline dynamics and $Z$ is known as the phase response curve (PRC), which describes the infinitesimal sensitivity of the phase to an external control input. One complete oscillation of the system corresponds to $\theta\in[0,2\pi)$. In the case of neural oscillators, $u$ represents an external current stimulus and $f$ is referred to as the instantaneous oscillation frequency in the absence of any external input, i.e., $u=0$.  As a convention, a neuron is said to spike or fire at time $T$ following a spike at time 0 if $\theta(t)$ evolves from $\theta(0)=0$ to $\theta(T)=2\pi$, i.e., spikes occur at $\theta=2n\pi$, where $n=0,1,2,\ldots$.  In the absence of any input $u(t)$ the neuron spikes periodically at its natural frequency, while by an appropriate choice of $u(t)$ the spiking time can be advanced or delayed in a desired manner. %Note that the conditions for validity and accuracy of the phase model have been determined \cite{efimov10}, and the reduction is accomplished through the well-studied process of phase coordinate transformation \cite{efimov09}, which is based on Floquet theory \cite{perko91,kelley04}.  The model is assumed valid for inputs $u(t)$ such that the solution $x(t,x_0,u)$ to (\ref{eq:sys1}) remains within a neighborhood of $\Gamma$.

In this article, we study various phase models characterized by different $f$ and $Z$ functions. In particular, we investigate the neural inputs that elicit desired spikes for an ensemble of isolated neurons with different natural dynamics, e.g., different oscillation frequencies. Fundamental questions on the controllability of these neuron systems and the design of optimal inputs that spike them arise naturally and will be discussed.

% =============with different natural frequencies==============================
% \subsection{Controllability of Uncoupled Neuron Assemblies}
\subsection{Controllability of Neuron Ensembles}
In this section, we analyze controllability properties of finite collections of neuron oscillators. We first consider the Theta neuron model (Type I neurons) which describes both superthreshold and subthreshold dynamics near a SNIPER (saddle-node bifurcation of a fixed
point on a periodic orbit) bifurcation \cite{Moehlis06, ermentrout96}.

%with angular frequency $\omega=2\sqrt{I_b}$.

% =================== Theta Neuron ==============
\subsubsection{Theta Neuron Model}
\label{sec:thetaneuron}
\begin{figure}[t]
\centering
\includegraphics[scale=0.2]{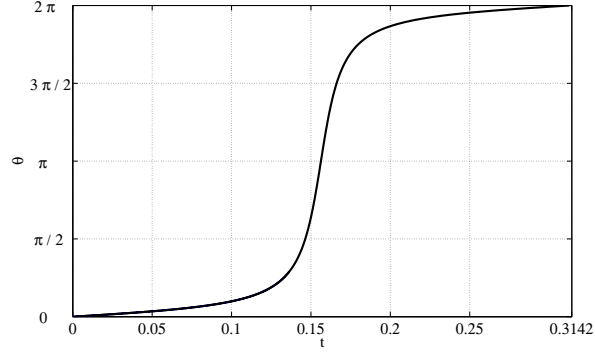}
\caption{The Free Evolution of a Theta Neuron. The baseline current $I=100$, and hence it spikes periodically with angular frequency $\w=20$ and period $T_0=\pi/10$.} \label{fig:thetaneuron_free}
\end{figure}

The Theta neuron model is characterized by the neuron baseline dynamics, $f(\theta)=(1+I)+(1-I)\cos\theta$, and the PRC, $Z(\theta)=1-\cos\theta$, namely,
% Consider the phase model as in \eqref{eq:phasemodel} with $f(\theta)=(1+I)+(1-I)\cos\theta$ and $Z(\theta)=1-\cos\theta$, namely,
\begin{eqnarray}
	\label{eq:thetaneuron1}
	\frac{d\t}{dt}=\big[(1+I)+(1-I)\cos\theta\big]+(1-\cos\theta)u(t),
\end{eqnarray}
where $I$ is the neuron baseline current. If $I>0$, then $f(\theta)>0$ for all $t\geq 0$. Therefore, in the absence of the input the neuron fires periodically since the free evolution of this neuron system, i.e., $\dot{\theta}=f(\theta)$, has a periodic orbit
\begin{eqnarray}
\label{eq:Ib+}
\theta(t)=2\tan^{-1}\left(\frac{\tan[\sqrt{I}(t+c)]}{\sqrt{I}}\right),\quad\text{for}\ \ I>0,
\end{eqnarray}
with the period $T_0=\pi/\sqrt{I}$ and hence the frequency $\w=2\sqrt{I}$, where $c$ is a constant depending on the initial condition. For example, if $\theta(0)=0$, then $c=0$. Fig. \ref{fig:thetaneuron_free} shows the free evolution of a Theta neuron with $I=100$. This neuron spikes periodically at $T_0=\pi/10$ with angular frequency $\w=2\pi/T=20$. When $I<0$, then the model is excitable, namely, spikes can occur with an appropriate input $u(t)$. However, no spikes occur without any input $u(t)$ as
$$\theta(t)=2\tan^{-1}\left(\frac{\tanh[\sqrt{-I}(t+c)]}{\sqrt{-I}}\right),\quad\text{for}\ \ I<0,$$
and there are two fixed points (one of which is stable) for $u(t)\equiv 0$.

%% =================== Figure 1 ==================

Now we consider spiking a finite collection of neurons with distinct natural oscillation frequencies and with positive baseline currents. This gives rise to a steering problem of the finite-dimensional single-input nonlinear control system, $\dot{\Theta}=f(\Theta)+Z(\Theta)u(t)$, where $\Theta\in\Omega\subset\mathbb{R}^n$, $f,Z:\Omega\rightarrow\mathbb{R}^n$, and $u\in U\subset\mathbb{R}$. In the vector form, this system appears as
\begin{eqnarray}
	\label{eq:neuronnetwork}
	\left[\begin{array}{c}\dot{\t}_1 \\ \dot{\t}_2 \\ \vdots \\
	\dot{\t}_n\end{array}\right]=\left[\begin{array}{c}\a_1+\b_1\cos\t_1 \\ \a_2+\b_2\cos\t_2 \\ \vdots \\
	\a_n+\b_n\cos\t_n\end{array}\right]+\left[\begin{array}{c}1-\cos\t_1 \\ 1-\cos\t_2 \\ \vdots \\
	1-\cos\t_n\end{array}\right]u(t),
\end{eqnarray}
in which $\a_i=1+I_i=1+\w^2_i/4$, $\b_i=1-I_i=1-\w^2_i/4$, and $I_i>0$ for all $i\in\mathcal{N}=\{1,2,\ldots,n\}$. Note that $\a_i+\b_i=2$ for all $i\in\mathcal{N}$. The ultimate proof of our understanding of neural systems is reflected in our ability to control them, hence a complete investigation of the controllability of oscillator populations is of fundamental importance. We now analyze controllability for the system as in \eqref{eq:neuronnetwork}, which determines whether spiking or synchronization of an oscillator ensemble by the use of an external stimulus is possible.

Because the free evolution of each neuron system $\theta_i$, $i\in\mathcal{N}$, in (\ref{eq:neuronnetwork}) is periodic as shown in (\ref{eq:Ib+}), the drift term $f$ causes no difficulty in analyzing controllability. The following theorem provides essential machinery for controllability analysis.

% =============== Theorem 1 ===============
\begin{thm}
	\label{thm:periodic}
	Consider the nonlinear control system
	\begin{eqnarray}
		\label{eq:nonsys}
		\dot{x}(t)=f(x(t))+u(t)\, g(x(t)),\quad x(0)=x_0.
	\end{eqnarray}
    Suppose that $f$ and $g$ are vector fields on a manifold $M$. Suppose that $\{f,g\}$ meet either of the conditions of Chow's theorem, and suppose that for each initial condition $x_0$ the solution of
	$$\dot{x}(t)=f(x(t))$$
    is periodic with a least period $T(x_0)<P\in\mathbb{R}^{+}$. Then the reachable set from $x_0$ for \eqref{eq:nonsys} is $\{\exp\{f,g\}_{LA}\}_G\, x_0$, where $\{f,g\}_{LA}$ denotes the Lie algebra generated by the vector fields $f$ and $g$, and $\{\exp\{f,g\}_{LA}\}_G$ is the smallest subgroup of the diffeomorphism group, ${\rm diff}(M)$, which contains $\exp t\eta$ for all $\eta\in\{f,g\}_{LA}$ \cite{Brockett76}.
\end{thm}
{\it Proof.}  See Appendix \ref{apd:Chow}. \hfill$\Box$

The underlying idea of this theorem for dealing with the drift is the utilization of periodic motions along the drift vector field, $f(x)$, to produce negative drift by forward evolutions for long enough time. More details about Theorem \ref{thm:periodic} can be found in Appendix \ref{apd:Chow}. Having this result, we are now able to investigate controllability of a neuron oscillator assembly.

% =============== Theorem 2 ===============
\begin{thm}
	\label{thm:thetaneuron}
	Consider the finite-dimensional single-input nonlinear control system
	\begin{eqnarray}
		\label{eq:thetaneuron2}
		\dot{\Theta}=f(\Theta)+Z(\Theta)u(t),
	\end{eqnarray}
	where $\Theta=(\t_1,\t_2,\ldots,\t_n)'\in\Omega\subset\mathbb{R}^n$ and the vector fields $f,Z:\Omega\rightarrow\mathbb{R}^n$ are defined by
	$$f(\Theta)=\left[\begin{array}{c}\a_1+\b_1\cos\t_1 \\ \a_2+\b_2\cos\t_2 \\ \vdots \\ \a_n+\b_n\cos\t_n\end{array}\right],\quad
Z(\Theta)=\left[\begin{array}{c}1-\cos\t_1 \\ 1-\cos\t_2 \\ \vdots \\ 1-\cos\t_n\end{array}\right],$$
	in which $\a_i=1+I_i$, $\b_i=1-I_i$, and $I_i>0$ for all $i\in\mathcal{N}=\{1,2,\ldots,n\}$. The system as in \eqref{eq:thetaneuron2} is controllable.
\end{thm}
%\begin{thm}
%The single-input nonlinear control system $\dot{\Theta}=f(\Theta)+g(\Theta)u(t)$, where $\Theta=(\t_1,\t_2,\ldots,\t_n)^{T}\in\Omega\subset\mathbb{R}^n$ and the vector fields $f$ and $g$ are defined as in \eqref{eq:neuronnetwork}, is controllable.
%\end{thm}
{\it Proof.} It is sufficient to consider the case where $I_i\neq I_j$ for $i\neq j$ and $i,j\in\mathcal{N}$, since otherwise they present the same neuron system. Because $I_i>0$ for all $i\in\mathcal{N}$, the free evolution, i.e., $u(t)=0$, of each $\t_i$ is periodic for every initial condition $\t_i(0)\in\mathbb{R}$, as shown in \eqref{eq:Ib+}, with the angular frequency $\w_i=2\sqrt{I_i}$ and the period $T_i=\pi/\sqrt{I_i}$. Therefore, the free evolution of $\Theta$ is periodic with a least period or is recurrent (see Remark \ref{rmk:period}). We may then apply Theorem \ref{thm:periodic} in computing the reachable set of this system. Let
$$ad_g h(\Theta)=[g,h](\Theta)$$
denote the Lie bracket of the vector fields $g$ and $h$, both defined on an open subset $\Omega$ of $\mathbb{R}^n$. Then, the recursive operation is denoted as
$$ad_g^k h(\Theta)=[g,ad_g^{k-1} h](\Theta)$$
for any $k\geq 1$, setting $ad_g^0 h(\Theta)=h(\Theta)$. The Lie brackets of $f$ and $Z$ include
\begin{align*}
	ad_f^{2k-1} Z &= \left[\begin{array}{c}(-1)^{k-1}2^k (\a_1-\b_1)^{k-1}\sin\t_1 \\ (-1)^{k-1}2^k (\a_2-\b_2)^{k-1}\sin\t_2 \\ \vdots \\ (-1)^{k-1}2^k (\a_n-\b_n)^{k-1}\sin\t_n\end{array}\right],\\
ad_f^{2k} Z &= \left[\begin{array}{c}(-1)^{k-1}2^k(\a_1-\b_1)^{k-1}(\a_1\cos\t_1+\b_1) \\ (-1)^{k-1}2^k(\a_2-\b_2)^{k-1}(\a_2\cos\t_2+\b_2) \\ \vdots \\ (-1)^{k-1}2^k(\a_n-\b_n)^{k-1}(\a_n\cos\t_n+\b_n)\end{array}\right],
\end{align*}
for $k\in\mathbb{Z}^{+}$, positive integers. Thus, $\{f,ad_f^m Z\}$, $m\in\mathbb{Z}^+$, spans $\mathbb{R}^n$ at all $\Theta\in\Omega$ since $\a_i-\b_i=\w_i^2/2$ and $\w_i$ are distinct for $i\in\mathcal{N}$. That is, every point in $\mathbb{R}^n$ can be reached from any initial condition $\Theta(0)\in\Omega$, hence the system \eqref{eq:thetaneuron2} is controllable. Note that if $\Theta=0$, $\{ad_f^{2k} Z\}$, $k\in\mathbb{Z}^+$, spans $\mathbb{R}^n$. $\hfill\Box$

% ================= Remark 1 ======================
\begin{rmk}
	\label{rmk:period}
	If there exist integers $m_i$ and $n_j$ such that the periods of neuron oscillators are related by $m_iT_i=n_j T_j$ for all $(i,j)$ pairs, $i,j\in\mathcal{N}$, then the free evolution of $\Theta$ is periodic with a least period. If, however, such a rational number relation does not hold between any two periods, e.g., $T_1=1$ and $T_2=\sqrt{2}$, it is easy to see that the free evolution of $\Theta$ is almost-periodic \cite{Levitan82} because the free evolution of each $\t_i$, $i\in\mathcal{N}$, is periodic. Hence, the recurrence of $f$ in \eqref{eq:thetaneuron2} together with the Lie algebra rank condition (LARC) described above guarantee the controllability.
	
	% If the frequency $\w_i$ of each system $\t_i$ as in \eqref{eq:thetaneuron2} is a rational number, then $\Theta$ is periodic with a least period equal to the least common multiple of $\{T_i=\frac{2\pi}{\w_i}\}$, $i\in\mathcal{N}$, where $T_i$ is the period of the system $\t_i$. The periodicity of $\Theta$ can be (approximately) preserved and Theorem \ref{thm:periodic} is still applied when $\w_i$ are irrational numbers, since these irrational $\w_i$ can be approximated by rational numbers with arbitrary precision.

% The periodicity of $\Theta$ in Theorem \ref{thm:thetaneuron} is preserved when $\w_i$ are irrational numbers, as what they usually are, since we can approximate any irrational numbers by rational numbers with arbitrary precision.
\end{rmk}

Controllability properties for other commonly-used phase models used to describe the dynamics of neuron or other, e.g., chemical, oscillators can be shown in the same fashion.

% ================== SNIPER PRC ===============
\subsubsection{SNIPER PRC}

\label{sec:sniper}
The SNIPER phase model is characterized by $f=\w$, the neuron's natural oscillation frequency, and the SNIPER PRC, $Z=z(1-\cos\t)$, where $z$ is a model-dependent constant \cite{Moehlis06}. In the absence of any external input, the neuron spikes periodically with the period $T=2\pi/\omega$. The SNIPER PRC is derived for neuron models near a SNIPER bifurcation which is found for Type I neurons \cite{ermentrout96} like the Hindmarsh-Rose model \cite{Rose89}. Note that the SNIPER PRC can be viewed as a special case of the Theta neuron PRC for the baseline current $I>0$. This can be seen through a bijective coordinate transformation $\t(\phi)=2\tan^{-1}[\sqrt{I_b}\tan(\phi/2-\pi/2)]+\pi$, $\phi\in[0,2\pi)$, applied to \eqref{eq:thetaneuron1}, which yields $\frac{d\phi}{dt}=\w+\frac{2}{\w}(1-\cos\phi)u(t)$, i.e., the SNIPER PRC with $z=2/\w$. The spiking property, namely, $\t(\phi=0)=0$ and $\t(\phi=2\pi)=2\pi$ is preserved under the transformation and so is the controllability as analyzed in Section \ref{sec:thetaneuron}.

More specifically, consider a finite collection of SNIPER neurons with $f(\Theta)=(\w_1,\w_2,\ldots,\w_n)'$ and $Z(\Theta)=(z_1(1-\cos\t_1),z_2(1-\cos\t_2),\ldots,z_n(1-\cos\t_n))'$, where conventionally, $z_i=2/\w_i$ for $i\in\mathcal{N}$. Similar Lie bracket computations as in the proof of Theorem \ref{thm:thetaneuron} result in, for $k=1,\ldots,n$,
\begin{align*}
	ad_f^{2k-1} Z &= \big((-1)^{k-1}z_1\w_1^{2k-1}\sin\t_1,\ldots,(-1)^{k-1}z_n\w_n^{2k-1}\sin\t_n\big)',\\
	ad_f^{2k} Z &= \big((-1)^{k-1}z_1\w_1^{2k}\cos\t_1,\ldots,(-1)^{k-1}z_n\w_n^{2k}\cos\t_n\big)',
\end{align*}
and thus ${\rm span}\{f,Z\}_{LA}=\mathbb{R}^n$, since $\w_i\neq\w_j$ for $i\neq j$. Therefore, the system of a network of SNIPER neurons is controllable.

% ================= Sinusoidal PRC ================
\subsubsection{Sinusoidal PRC}
\label{sec:sinusoidal}
In this case, we consider $f(\Theta)=(\w_1,\w_2,\ldots,\w_n)'$ and $Z(\Theta)=(z_1\sin\t_1,$ $z_2\sin\t_2,\ldots,z_n\sin\t_n)'$, where $\w_i>0$ and $z_i=2/\w_i$ for $i=1,\ldots,n$. This type of PRC's with both positive and negative regions can be obtained by periodic orbits near the super critical Hopf bifurcation\cite{Brown04}. This type of bifurcation occurs for Type II neuron models like Fitzhugh-Nagumo model \cite{Keener98}. Controllability of a network of Sinusoidal neurons can be shown by the same construction, from which
\begin{align*}
	ad_f^{2k-1} Z &= \big((-1)^{k-1}z_1\w_1^{2k-1}\cos\t_1,\ldots,(-1)^{k-1}z_n\w_n^{2k-1}\cos\t_n\big)',\\
	ad_f^{2k} Z &= \big((-1)^{k}z_1\w_1^{2k}\sin\t_1,\ldots,(-1)^{k}z_n\w_n^{2k}\sin\t_n\big)',
\end{align*}
and then ${\rm span}\{f,Z\}_{LA}=\mathbb{R}^n$ for $\w_i\neq\w_j$, $i\neq j$. Therefore the system is controllable.

%%%%%%%%%%%%%%%%%%%%%%%%%%%%%%%%%%%%%%%%%%%%%%%%%%%%%%%%%%%%%%%%%%%%%
\section{Optimal Control of Spiking Neurons}
\label{sec:oc}
The controllability addressed above guarantees the existence of an input that drives an ensemble of oscillators between any desired phase configurations. Practical applications demand minimum-power or time-optimal controls that form certain synchronization patterns for a population of oscillators, which gives rise to an optimal steering problem,

% The controllability addressed above guarantees the existence of control inputs or external stimuli that lead to simultaneous spikes of a neuron ensemble. Practical applications, such as deep brain stimulations and cardiac pacemakers, make it desirable to find optimal controls, e.g., minimum-power stimuli, that synchronize a neuron ensemble. This then forms an optimal steering problem,

\begin{align}
	% \min\quad & \int_0^{t_1}|u(t)|^2 dt\nonumber\\
	\min\quad & J=\varphi(T,\Theta(T))+\int_0^{T}\mathcal{L}(\Theta(t),u(t)) dt\nonumber\\
	{\rm s.t.}\quad & \dot{\Theta}(t)=f(\Theta)+Z(\Theta)u(t) \label{eq:steering}\\
	& \Theta(0)=\Theta_0,\ \Theta(T)=\Theta_d \nonumber\\ %2\pi \text{P}.\nonumber\\
	& |u(t)|\leq M,\quad \forall\, t\in [0,T],\nonumber
\end{align}
where $\Theta\in\bR^n$, $u\in\bR$; $\varphi:\bR\times\bR^n\rightarrow\bR$, denoting the terminal cost, $\mathcal{L}:\bR^n\times\bR\rightarrow\bR$, denoting the running cost, and $f,Z:\bR^n\rightarrow\bR^n$ are Lipschitz continuous (over the respective domains) with respect to their arguments. For spiking a neuronal population, for example, the goal is to drive the system from the initial state, $\Theta_0=\textbf{0}$, to a final state $\Theta_d=(2m_1\pi,2m_2\pi,\ldots,2m_n\pi)'$, where $m_i\in\mathbb{Z}^+$, $i=1,\ldots, n$. Steering problems of this kind have been well studied, for example, in the context of nonholonomic motion planning and sub-Riemannian geodesic problems \cite{Sastry92, Montgomery89}. This class of optimal control problems in principle can be approached by the maximum principle, however, in most cases they are analytically intractable especially when the system is of high dimension, e.g., greater than three, and when the control is bounded, i.e., $M<\infty$. In the following, we present analytical optimal controls for single- and two-neuron systems and, furthermore, develop a robust computational method for solving challenging optimal control problems of steering a neuron ensemble. Our numerical method is based on pseudospectral approximations which can be easily extended to consider any topologies of neural networks, e.g., arbitrary frequency distributions and coupling strengths between neurons, with various types of cost functional.

% ================== Single Neuron ===============
\subsection{Minimum-Power Control of a Single Neuron Oscillator}
Designing minimum-power stimuli to elicit spikes of neuron oscillators is of clinical importance, such as deep brain stimulation, used for a variety of neurological disorders including Parkinson's disease, essential tremor, and Dystonia, and neurological implants of cardiac pacemakers, where mild stimulations and low energy consumption are required \cite{Benabid91, Marks05}. Optimal controls for spiking a single neuron oscillator can be derived using the maximum principle. In order to illustrate the idea, we consider spiking a Theta neuron, described in \eqref{eq:thetaneuron1}, with minimum power. In this case, the cost functional is $J=\int_0^{T}u^2(t)dt$, and the initial and target states are 0 and $2\pi$, respectively. We first examine the case when the control is unbounded.

% --------- unbounded case --------------
The control Hamiltonian of this optimal control problem is defined by $H=u^2+\lambda(\alpha+\beta\cos\theta+u-u\cos\theta)$,
% \begin{equation}
%     \label{eq:hamiltonian}
% 	H=u^2+\lambda(\alpha+\beta\cos\theta+u-u\cos\theta),
% \end{equation}
where $\lambda$ is the Lagrange multiplier. The necessary conditions for optimality yield $\dot{\lambda}=-\frac{\partial H}{\partial \theta}=\lambda(\beta-u)\sin\theta$,
% \begin{align}
% 	\label{eq:lambda_dot}
% 	\dot{\lambda} &= -\frac{\partial H}{\partial \theta}=\lambda(\beta-u)\sin\theta,\\
% 	\frac{\partial H}{\partial u} &= 2u+\lambda(1-\cos\theta)=0.\nonumber
% \end{align}
and $u=-\frac{1}{2}\lambda (1-\cos\theta)$ by $\frac{\partial H}{\partial u}=0$. With these conditions, the optimal control problem is then transformed to a boundary value problem, which characterizes the optimal trajectories of $\theta(t)$ and $\lambda(t)$. We then can derive the optimal feedback law for spiking a Theta neuron at the specified time $T$ by solving the resulting boundary value problem,
\begin{equation}
    \label{eq:I*}
    u^*(\t)= \frac{-(\alpha+\beta\cos\theta)+\sqrt{(\alpha+\beta\cos\theta)^2-2\lambda_0(1-\cos\theta)^2}}{1-\cos\t},
\end{equation}
where $\lambda_0=\lambda(0)$, which can be obtained according to
\begin{equation}
	\label{eq:T}
	T=\int_0^{2\pi}{\frac{1}{\sqrt{(\alpha+\beta\cos\theta)^2-2\lambda_0(1-\cos\theta)}}}d\theta.
	% T=\int_0^{2\pi}{\frac{1}{\scriptstyle{\sqrt{(\alpha+\beta\cos\theta)^2-2\lambda_0(1-\cos\theta)}}}}d\theta.
\end{equation}
More details about the derivations can be found in Appendix \ref{apd:thetaneuron_unbounded}.

% --------- bounded case --------------
Now consider the case when the control amplitude is limited, namely, $|u(t)|\leq M$, $\forall\, t\in[0,T]$. If the unbounded minimum-power control as in \eqref{eq:I*} satisfies $|u^*(\t)|\leq M$ for all $t\in[0,T]$, then the amplitude constraint is inactive and obviously the optimal control is given by \eqref{eq:I*} and \eqref{eq:T}. However, if $|u^*(\t)|> M$ for some $\theta\in[0,2\pi]$, then the optimal control $u^*_{M}$ is characterized by switching between $u^*(\t)$ and the bound $M$ (see Appendix \ref{apd:thetaneuron_bounded}),
\begin{equation}
	\label{eq:bounded_control}
    u^*_{M}(\theta)=\left\{\begin{array}{ll}-M, & u^*(\theta)< -M \\ u^*(\theta), & -M\leq u^*(\theta)\leq M\\ M, & u(\theta)^*> M,
    \end{array}\right.
\end{equation}
where the parameter $\lambda_0$ for $u^*(\t)$ in \eqref{eq:I*} is calculated according to the desired spiking time $T$ by
\begin{equation}
\label{eq:bounded_control_T}
T=\int_0^{2\pi}\frac{1}{\a+\b\cos\t+(1-\cos\t)u^*_M}d\theta.
\end{equation}

The detailed derivation of the control $u_M^*$ is given in Appendix \ref{apd:thetaneuron_bounded}. Fig. \ref{fig:single_theta_neuron} illustrates the optimal controls and the corresponding trajectories for spiking a Theta neuron with natural oscillation frequency $\w=1$, i.e., $I=0.25$, $\alpha=1.25$, and $\beta=0.75$, at various spiking times that are smaller $(T=3,4)$ and greater $(T=8)$ than the natural spiking time $T_0=2\pi$ with the control amplitude bound $M=1$. Because the unconstrained minimum-power controls for the cases $T=4<T_0$ and $T=8>T_0$, calculated according to \eqref{eq:I*}, satisfy $|u^*(\t)|<M$, there are no switchings in these two cases.

 % and the natural oscillation frequency $\w=1$, which is equivalent to $I=0.25$, $\alpha=1.25$, and $\beta=0.75$.

% The detailed derivation of the control $u_M^*$ is given in Appendix \ref{apd:thetaneuron_bounded}. Fig. \ref{fig:single_theta_neuron} illustrates the optimal controls and the corresponding trajectories for spiking a Theta neuron at various spiking times $T=3,4,8,10$ with the control amplitude bound $M=1$ and the natural oscillation frequency $\w=1$, which is equivalent to $I=0.25$, $\alpha=1.25$, and $\beta=0.75$.
% % In this article, we apply techniques from optimal control theory to derive minimum-power controls that spike a neuron at desired time instants. We consider bounded control amplitude and fully characterize the range of feasible spiking times determined by the bound. In particular, our optimal control strategies are general so that the bound can be chosen sufficiently small within the range that the PRC is valid. The design of such minimum-power stimuli to elicit spikes of neuron oscillators is also of clinical importance, notably in deep brain stimulation therapy for Parkinson's disease and essential tremor \cite{Benabid10}, where mild stimulations are required. In addition, interest of reducing the energy consumption in neurological implants such as cardiac pacemakers makes such optimal designs imperative.

% ================= Figure 2 ===============
\begin{figure}[t]
	\subfigure[]{\includegraphics[scale=0.19]{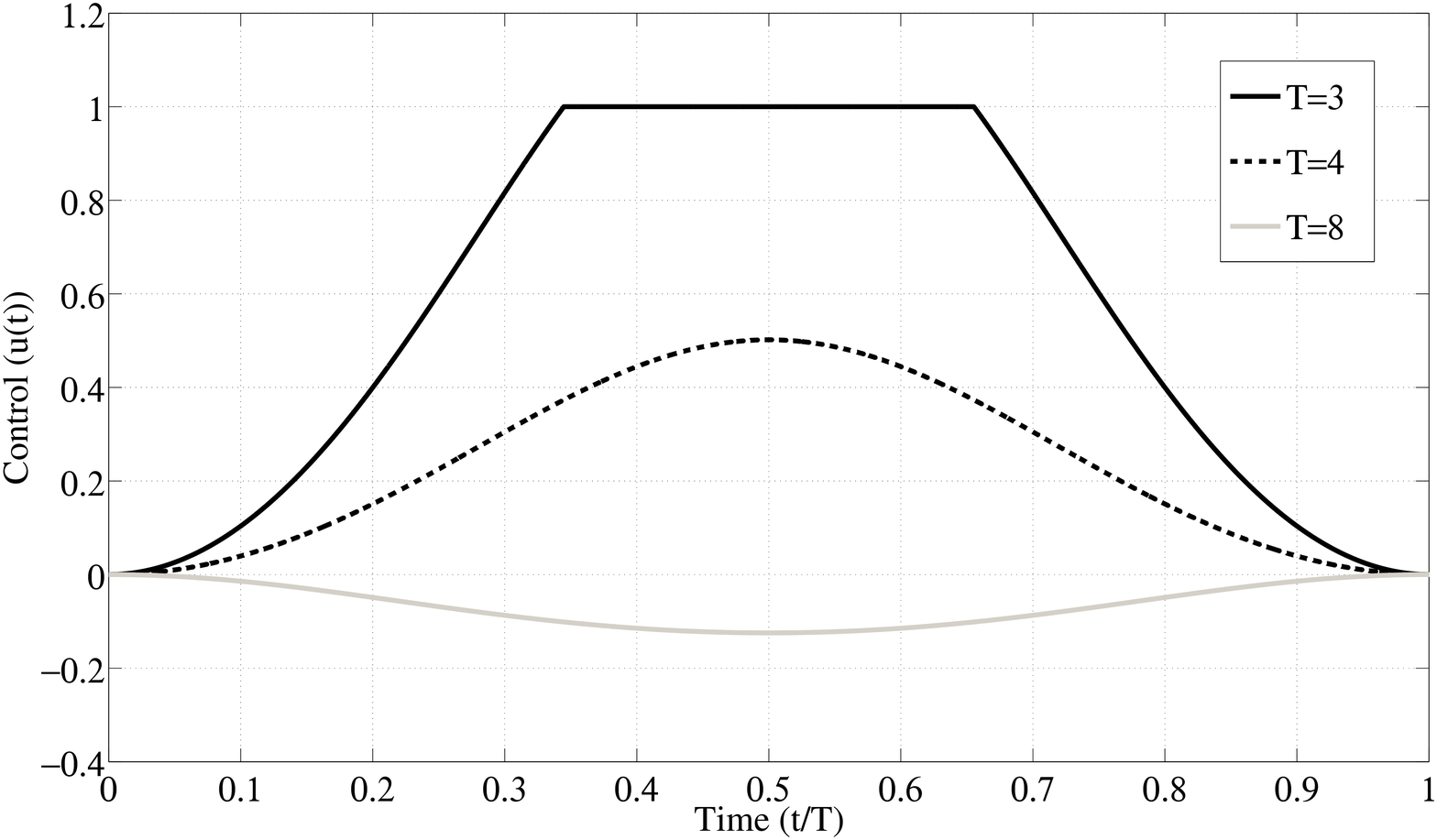}\label{fig:control_vs_time_theta}}
	\subfigure[]{\includegraphics[scale=0.19]{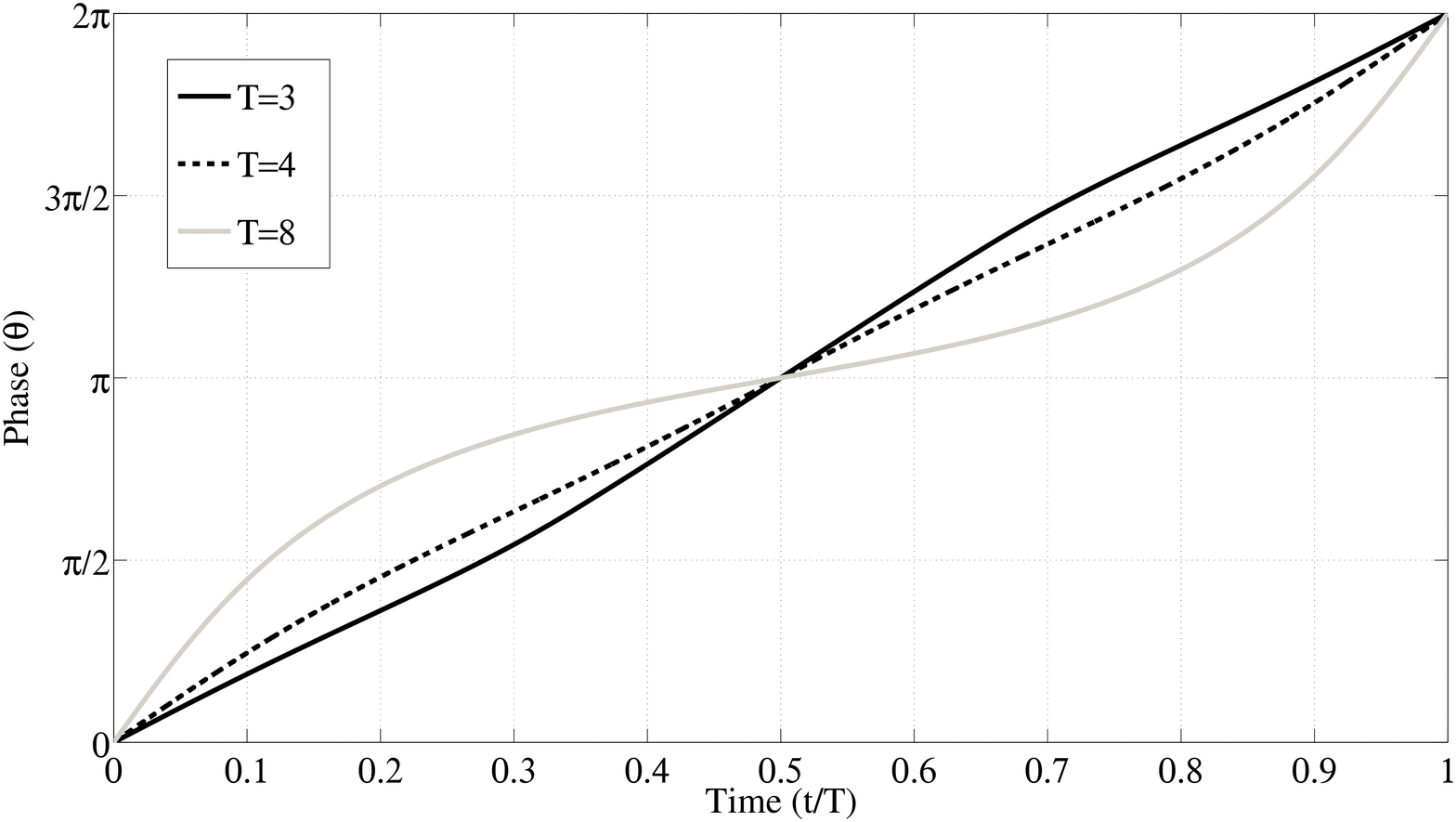}\label{fig:trajectory_theta}}
	\caption{\subref{fig:control_vs_time_theta} Minimum-power controls, $u^*_M(t)$, for spiking a Theta neuron with $\w=1$ at various spiking times that are smaller $(T=3,4)$ and greater $(T=8)$ than the natural spiking time $T_0=2\pi$ subject to $M=1$. \subref{fig:trajectory_theta} The resulting optimal phase trajectories following $u^*_M(t)$.}
	\label{fig:single_theta_neuron}
	% \caption{\subref{fig:control_vs_time_theta} Minimum-power controls, $u^*_M(t)$, for spiking a Theta neuron at various spiking times $T=3,4,8,10$ ms with $\omega=1$ and $M=1$. \subref{fig:trajectory_theta} The resulting optimal phase trajectories following $u^*(t)$.}
	% 	\label{fig:single_theta_neuron}
\end{figure}

% ================== Two Neurons ===============
\subsection{Time-Optimal Control of Two Neuron Oscillators}
Spiking a neuron in minimum time, subject to a given control amplitude, can be solved in a straightforward manner. Consider the phase model of a single neuron as in \eqref{eq:phasemodel}, it is easy to see that for a given control bound $M>0$, the minimum spiking time is achieved by the bang-bang control
\begin{equation}
	\label{eq:u_Tmin}
    u^*_t=\left\{\begin{array}{ll} M, & Z(\theta)\geq 0 \\ -M, & Z(\theta)<0, 
    \end{array}\right.
\end{equation}
which keeps the phase velocity, $\dot{\t}$, at its maximum. The minimum spiking time with respect to the control bound $M$, denoted by $T_{min}^M$, is then given by
\begin{align}
	T_{min}^M=&\int\limits_{\theta\in\mathcal{A}}\frac{1}{f(\theta)+Z(\theta)M}d\theta+\int\limits_{\theta\in\mathcal{B}} \frac{1}{f(\theta)-Z(\theta)M}d\theta,
\label{TminM}
\end{align}
where the sets $\mathcal{A}$ and $\mathcal{B}$ are defined as
\begin{align*}
	\mathcal{A}=&\left\{\theta|\ Z(\theta)\geq 0, 0\leq \theta \leq 2\pi\right\},\\
	\mathcal{B}=&\left\{\theta|\ Z(\theta)< 0, 0\leq \theta \leq 2\pi\right\}.
\end{align*}
Time-optimal control of spiking two neurons is more involved, which can be formulated as in \eqref{eq:steering} with the cost functional $J=\int_0^{T} 1 dt$ and with 
\begin{eqnarray}
	\label{eq:twoneuron}
	\dot{\Theta}(t)=f(\Theta)+Z(\Theta)u(t),
\end{eqnarray}
where
\begin{eqnarray}
	\label{eq:twoneuron1}
	f=\left[\begin{array}{c}f_1 \\ f_2\end{array}\right]=\left[\begin{array}{c}\a_1+\b_1\cos\t_1 \\ \a_2+\b_2\cos\t_2\end{array}\right],
	\quad Z=\left[\begin{array}{c} Z_1 \\ Z_2\end{array}\right]=\left[\begin{array}{c} 1-\cos\t_1 \\ 1-\cos\t_2\end{array}\right].
\end{eqnarray}
Our objective is to drive the two-neuron system from the initial state $\Theta_0=(0, 0)'$ to the desired final state $\Theta_d=(2m_1\pi,2m_2\pi)'$ with minimum time, where $m_1,m_2\in\mathbb{Z}^{+}$. 
%Time-optimal control of spiking two neurons is more involved, which can be formulated as in \eqref{eq:steering} with the cost functional $J=\int_0^{T} 1 dt$. Let's consider spiking two theta neurons simultaneously in minimum-time, where $f=(\a_1+\b_1\cos\t_1,\a_2+\b_2\cos\t_2)'$,  $Z=(1-\cos\t_1,1-\cos\t_2)'$. Our objective is to drive the two-neuron system from the initial state $(0, 0)'$ to the final state $(2\pi n,2\pi m)'$ with minimum time, where $n,m\in\mathbb{Z}^{+}$. 
The Hamiltonian for this optimal control problem is given by
\begin{eqnarray}
	\label{eq:hamiltonian2}
	H=\lambda_0+\<\lambda,f+Zu\>,
\end{eqnarray}
where $\lambda_0\in\mathbb{R}$ and $\lambda\in\mathbb{R}^2$ are the multipliers that correspond to the Lagrangian and the system dynamics, respectively, and $\<\, , \,\>$ denotes a scalar product in the Euclidean space $\mathbb{E}^2$.

% ======== Proposition 1 ===========
\begin{prop}
	The minimum-time control that spikes two Theta neurons simultaneously is bang-bang.
\end{prop}
{\it Proof.} %We prove this by examining two cases. 
The Hamiltonian in \eqref{eq:hamiltonian2} is minimized by the control,
\begin{equation}
	\label{eq:opt_control_2_theta}
	u(t)=\left\{\begin{array}{lll} M &\text{for} & \phi(t)< 0, \\ -M &\text{for} & \phi(t)>0, \end{array} \right.
\end{equation}
where $\phi$ is the switching function defined by $\phi=\<\lambda,Z\>$. If there exists no non-zero time interval over which $\phi\equiv 0$, then the optimal control is given by the bang-bang form as in \eqref{eq:opt_control_2_theta}, where the control switchings are defined at $\phi=0$. We show by contradiction that maintaining $\phi=0$ is not possible for any non-zero time interval. Suppose that $\phi(t)=0$ for some non-zero time interval, $t\in[\tau_1,\tau_2]$, then we have
\begin{align}
	\label{phi}
	\phi =& \<\lambda,Z\>=0,\\
	\label{phidot}
	\dot{\phi} =& \<\lambda,\,[f,Z]\,\>=0,
\end{align}
where $[f,Z]$ denotes the Lie bracket of the vector fields $f$ and $Z$.
% \noindent{\it Case I: $f$ and $Z$ are linearly independent}
According to  \eqref{phi} and \eqref{phidot}, $\lambda$ is perpendicular to both vectors $Z$ and $[f,Z]$, where
\[
[f,Z] = \frac{\partial Z}{\partial\t} f-\frac{\partial f}{\partial\t}Z=\left[\begin{array}{c}2\sin\theta_1 \\ 2\sin\theta_2 \\ \end{array}\right].
\]
Since $\lambda\neq0$ by the non-triviality condition of the maximum principle, $Z$ and $[f,Z]$ are linearly dependent on $t\in[\tau_1,\tau_2]$. One can easily show that these two vectors are linearly dependent either when $\t_1=2n\pi$ and $\t_2\in\mathbb{R}$, $\t_1\in\mathbb{R}$ and $\t_2=2n\pi$, or $\t_1=\t_2+2n\pi$ and $\t_2\in\mathbb{R}$, where $n\in \mathbb{Z}$. These three families of lines represent the possible paths in the state-space where $\phi$ can be vanished for some non-trivial time-interval. Now we show that these are not feasible phase trajectories that can be generated by a control. Suppose that $(\t_1(\tau),\t_2(\tau))=(2n\pi,\a)$ for some $\tau>0$ and for some $n\in\mathbb{Z}$, where $\a\in\mathbb{R}$. We then have $\dot{\t_1}(\tau)=2\neq 0$, irrespective of any control input. Hence, the system is immediately deviated from the line $\t_1=2n\pi$. The same reasoning can be used for showing the case of $\t_2=2n\pi$.

Similarly, if $(\t_1(\tau),\t_2(\tau))=(\a+2n\pi,\a)$ for some $\tau>0$ and for some $n\in\mathbb{Z}$, in order for the system to remain on the line $(\t_1(t),\t_2(t))=(\t_2+2n\pi,\t_2)$, it requires that $\dot{\t_1}(t)=\dot{\t_2}(t)$ for $t>\tau$. However, this occurs only when $\t_1=2m\pi$ and $\t_2=2(n+m)\pi$, where $m\in\mathbb{Z}$, since $\dot{\t_1}-\dot{\t_2}=(I_{1}-I_{2})(1-\cos\t_1)$. Furthermore, staying on these points is impossible with any control inputs since for $\t_1(\tau)=2m\pi$ and $\t_2(\tau)=2(n+m)\pi$, the phase velocities are $\dot{\t_1}(\tau)=\dot{\t_2}(\tau)=2$, which immediately forces the system to be away from these points. Therefore, the system cannot be driven along the path $(\t_2+2n\pi,\t_2)$. This analysis concludes that $\phi=0$ and $\dot{\phi}=0$ do not hold simultaneously over a non-trivial time interval. \hfill$\Box$

Now, we construct the bang-bang structure for time-optimal control of this two-neuron system and, without loss of generality, let $\lambda_0=1$.

\begin{defn}
	We denote the vector fields corresponding to the constant bang controls $u(t)\equiv -M$ and $u(t)\equiv M$ by $X = f - MZ$ and $Y = f + MZ$, respectively, and call the respective trajectories corresponding to them as $X$- and $Y$- trajectories. A concatenation of an $X$-trajectory followed by a $Y$-trajectory is denoted by $XY$, while the concatenation in the reverse order is denoted by $YX$. 
\end{defn}

Due to the bang-bang nature of the time-optimal control for this system, it is sufficient for us to calculate the time between consecutive switches, and then the first switching time can be determined by the end point constraint. The inter-switching time can be calculated following the procedure described in \cite{Sussmann82,Sussmann87,Li_SICON11}. 

Let $p$ and $q$ be consecutive switching points, and let $pq$ be a $Y$-trajectory. Without loss of generality, we assume that this trajectory passes through $p$ at time $0$ and is at $q$ at time $\tau$. Since $p$ and $q$ are switching points, the corresponding multipliers vanish against the control vector field $Z$ at those points, i.e., 
\begin{eqnarray}
	\label{eq:lambdaZ_1}
	\<\lambda(0),Z(p)\>=\<\lambda(\tau), Z(q)\>=0.
\end{eqnarray}
Assuming that the coordinate of $p=(\theta_1,\theta_2)'$, our goal is to calculate the switching time, $\tau$, in terms of $\theta_1$ and $\theta_2$. In order to achieve this, we need to compute what the relation $\<\lambda(\tau), Z(q)\> = 0$ implies at time $0$. This can be obtained by moving the vector $Z(q)$ along the $Y$-trajectory backward from $q$ to $p$ through the pushforward of the solution $\omega(t)$ of the variational equation along the $Y$-trajectory with the terminal condition $\omega(\tau)=Z(q)$ at time $\tau$. We denote by $e^{tY}(p)$ the value of the $Y$-trajectory at time $t$ that starts at the point $p$ at time $0$ and by $(e^{-tY})_{\ast}$ the backward evolution under the variational equation. Then we have
$$\w(0)=(e^{-\tau Y})_{\ast}\,\w(\tau)=(e^{-\tau Y})_{\ast}\,Z(q)=(e^{-\tau Y})_{\ast}\,Z(e^{\tau Y}(p))=(e^{-\tau Y})_{\ast}\, Z\, e^{\tau Y}(p).$$
% \begin{eqnarray}
% 	\w(0)=(e^{-\tau Y})_{\ast}\,\w(\tau)=(e^{-\tau Y})_{\ast}\,Z(q)=(e^{-\tau Y})_{\ast}\,Z(e^{\tau Y}(p))=(e^{-\tau Y})_{\ast}\, Z\, e^{\tau Y}(p).
% \end{eqnarray}
Since the ``adjoint equation'' of the maximum principle is precisely the adjoint equation to the variational equation, it follows that
the function $t\mapsto\langle\lambda(t),\w(t)\rangle$ is constant along the $Y$-trajectory. Therefore, $\langle\lambda(\tau),Z(q)\rangle=0$ also implies that
\begin{eqnarray}
	\label{eq:lambdaZ_2}
	\<\lambda(0),\w(0)\>=\<\lambda(0),(e^{-\tau Y})_{\ast}\,Z\, e^{\tau Y}(p)\>=0.
\end{eqnarray}
Since $\lambda(0)\neq 0$, we know from \eqref{eq:lambdaZ_1} and \eqref{eq:lambdaZ_2} that the two vectors $Z(p)$ and $(e^{-\tau Y})_{\ast}\,Z\, e^{\tau Y}(p)$ are linearly dependent. It follows that 
\begin{equation}
	\label{eq:gamma}
	\gamma Z(p)=(e^{-\tau Y})_{\ast}\,Z\, e^{\tau Y}(p),
\end{equation}
where $\gamma$ is a constant. We make use of a well-known Campbell-Baker-Hausdorff formula \cite{Isidori95} to expand $(e^{-\tau Y})_{\ast}\,Z\, e^{\tau Y}(p)$, that is,
\begin{equation*}
	(e^{-\tau Y})_{\ast}\,Z\, e^{\tau Y}(p)=e^{\tau ad_Y}(Z)=\sum_{n=0}^{\infty}\frac{\tau^n}{n!}ad^n_YZ.
\end{equation*}
% \begin{equation*}
% 	pe^{\tau Y}Ze^{-\tau Y}=e^{\tau ad_Y}(Z)
% \end{equation*}
A straightforward computation of Lie brackets gives
\begin{align*}
	ad_Y Z &= [Y,Z]=[f+MZ,Z]=[f,Z]=2\left[\begin{array}{l} \sin\theta_1 \\ \sin\theta_2 \end{array}\right],\\
	ad_Y^2Z &= [Y,[Y,Z]]=2(f-AZ),
\end{align*}
where $A=\text{diag} \left\{ 2(\alpha_1-2+M),  2(\alpha_2-2+M)\right\}$, and furthermore
\begin{eqnarray*}
ad^{2n+1}_YZ=(-1)^n2^n(A+MI)^n[f,Z],\\
ad^{2n+2}_YZ=(-1)^n2^{n+1}(A+MI)^n(f-AZ).
\end{eqnarray*}
Consequently, we have
\begin{align*}
	\scriptstyle{e^{\tau ad_YZ}=Z+\sum_{n=0}^{\infty}\frac{\tau^{2n+1}}{(2n+1)!}(-1)^n2^n(A+MI)^n[f,Z] +\sum_{n=0}^{\infty}\frac{\tau^{2n+2}}{(2n+2)!}(-1)^{n}2^{n+1}(A+MI)^n(f-AZ),}
\end{align*}
which is further simplified to
\begin{align*}
	\scriptstyle{e^{\tau ad_YZ}}=
\left[\begin{array}{c}\scriptstyle{\frac{\alpha_1+M-(M-\beta_1)\cos\theta_1}{2(\alpha_1-1+M)}+\frac{M-\beta_1-(\alpha_1+M)\cos\theta_1}{2(\alpha_1-1+M)}\cos(2\tau\sqrt{\scriptstyle{\alpha_1-1+M}})+\frac{\sin\theta_1}{\sqrt{\alpha_1-1+M}}\sin(2\tau\sqrt{\scriptstyle{\alpha_1-1+M}})}\\
	\scriptstyle{\frac{\alpha_2+M-(M-\beta_2)\cos\theta_2}{2(\alpha_2-1+M)}+\frac{M-\beta_2-(\alpha_2+M)\cos\theta_2}{2(\alpha_2-1+M)}\cos(2\tau\sqrt{\scriptstyle{\alpha_2-1+M}})+\frac{\sin\theta_2}{\sqrt{\alpha_2-1+M}}\sin(2\tau\sqrt{\scriptstyle{\alpha_2-1+M}})}\end{array}\right].
\end{align*}
This together with \eqref{eq:gamma} yields
\begin{align}
	&\scriptstyle{(1-\cos\theta_2)\left[\frac{\alpha_1+M-(M-\beta_1)\cos\theta_1}{2(\alpha_1-1+M)}+\frac{M-\beta_1-(\alpha_1+M)\cos\theta_1}{2(\alpha_1-1+M)}\cos(2\tau\sqrt{\alpha_1-1+M})+\frac{\sin\theta_1}{\sqrt{\alpha_1-1+M}}\sin(2\tau\sqrt{\alpha_1-1+M})\right]}\nonumber\\
	\label{eq:interswitching}
	&\scriptstyle{=(1-\cos\theta_1)\left[\frac{\alpha_2+M-(M-\beta_2)\cos\theta_2}{2(\alpha_2-1+M)}+\frac{M-\beta_2-(\alpha_2+M)\cos\theta_2}{2(\alpha_2-1+M)}\cos(2\tau\sqrt{\alpha_2-1+M})+\frac{\sin\theta_2}{\sqrt{\alpha_2-1+M}}\sin(2\tau\sqrt{\alpha_2-1+M})\right].}
%(\alpha-1)(1-\cos\theta_2)=\frac{(\alpha_2+M)\cos\theta_2+\beta-M}{M(2\alpha_2-1)}(1-\cos(\tau\sqrt{2M(\alpha_2-1)}))\nonumber\\
%+\frac{2\sin\theta_2}{\sqrt{2M(2\alpha_2-1)}}\sin(\tau\sqrt{2M(2\alpha_2-1)}). \label{eq:tau2}
\end{align}
This equation characterizes the inter-switching along the $Y$-trajectory, that is, the next switching time $\tau$ can be calculated given the system starting with $(\t_1,\t_2)$ evolving along the $Y$-trajectory. Similarly, the inter-switching along the $X$-trajectory can be calculated by substituting $M$ with $-M$ in \eqref{eq:interswitching}.

Note that the solution to \eqref{eq:interswitching} is not unique, and some of the solutions may not be optimal, which can be discarded in a systematic way. The idea is to identify those possible switching points calculated from \eqref{eq:interswitching} with $\phi=0$ that also having the appropriate sign for $\dot{\phi}$. We focus on the case where $f$ and $Z$ are linearly independent, since the case for those being linearly dependent restricts the state space to be the curve
% since the case for those being linearly dependent is very restricted, characterized by the curve
\begin{eqnarray*}
	% \label{eq:abnormal_curve}
	(\alpha_1+\beta_1\cos\theta_1)(1-\cos\theta_2)=(\alpha_2+\beta_2\cos\theta_2)(1-\cos\theta_1).
\end{eqnarray*}
If $f$ and $Z$ are linearly independent, then $[f,Z]$ can be written as $[f,Z]=k_1 f+k_2 Z$, where
% We focus on the case where $f$ and $Z$ are linearly independent, since the case for those are linearly dependent is very restricted as characterized by \eqref{eq:abnormal_curve}. If $f$ and $Z$ are linearly independent, then $[f,Z]$ can be written as $[f,Z]=k_1 f+k_2 Z$, where
\begin{align*}
	k_1=&\frac{2\sin\t_1(1-\cos\t_2)-2\sin\t_2(1-\cos\t_1)}{(\alpha_1+\beta_1\cos\t_1)(1-\cos\t_2)-(\alpha_2+\beta_2\cos\t_2)(1-\cos\t_1)},\\
	k_2=&\frac{2\sin\t_1(\alpha_1+\beta_1\cos\t_1)-2\sin\t_2(\alpha_2+\beta_2\cos\t_2)}{(\alpha_1+\beta_1\cos\t_1)(1-\cos\t_2)-(\alpha_2+\beta_2\cos\t_2)(1-\cos\t_1)}.
\end{align*}
As a result, we can write $\dot{\phi}=\<\lambda,[f,Z]\>=k_1\<\lambda,f\>+k_2\<\lambda,Z\>$. Since we know that at switching points $\phi=\<\lambda,Z\>=0$, the Hamiltonian, as in \eqref{eq:hamiltonian2}, $H=0$ and the choice of $\lambda_0=1$ makes $\<\lambda,f\>=-1$. Therefore, at these points, we have $\dot{\phi}=-k_1$, and the type of switching can be determined according to the sign of the function $k_1$. If $k_1>0$, then it is an $X$ to $Y$ switch since $\dot{\phi}<0$ and hence $\phi$ changes its sign from positive to negative passing through the switching point, which corresponds to switch the control from $-M$ to $M$ as in \eqref{eq:opt_control_2_theta}. Similarly, if $k_1<0$, then it is a $Y$ to $X$ switch. Therefore the next switching time will be the minimum non-zero solution to the equation \eqref{eq:interswitching} that satisfy the above given rule. For example, suppose that the system is following a $Y$-trajectory starting with a switching point $p_i=(\t^i_{1},\t^i_{2})'$. The possible inter-switching times $\{\tau_{i,j}\}$, $j=1,\ldots,n$, with $\tau_{i,1}<\tau_{i,2}<\ldots<\tau_{i,n}$ can then be calculated according to \eqref{eq:interswitching} based on $p_i$. Thus, the next switching point is $p_{r}=(\t^{r}_1,\t^{r}_2)'=e^{\tau_{i,r} Y}(p_i)$, $\tau_{i,r}=\min\{\tau_{i,1},\ldots,\tau_{i,n}\}$, such that $k_1(\t^r_1,\t^r_2)<0$, which corresponds to an $Y$ to $X$ switch.
% The possible inter-switching times $\{\tau_j\}$, $j=1,\ldots,n$, with $\tau_{1}<\tau_{2}<\ldots<\tau_{n}$ can then be calculated according to \eqref{eq:interswitching} based on $p_i$. Thus, the next switching point is $p_{r}=(\t^{r}_1,\t^{r}_2)'=e^{\tau_{r} Y}(p_i)$, $\tau_r=\min\{\tau_1,\ldots,\tau_n\}$, such that $k_1(\t^r_1,\t^r_2)<0$, which corresponds to an $Y$ to $X$ switch.

% For example, suppose that the system starts with a switching point $p_0=(\t^0_{1},\t^0_{2})'$ and follow a $Y$-trajectory. The possible inter-switching times $\{\tau_i\}$, $i=1,\ldots,n$, with $\tau_1<\tau_2<\ldots<\tau_n$ can then be calculated according to \eqref{eq:interswitching} based on $p_0$. Thus, the next switching point is $p_i=(\t^i_1,\t^i_2)'=e^{\tau_i Y}(p_0)$ such that $k_1(\t^i_1,\t^i_2)<0$, which corresponds to an $Y$ to $X$ switch.

Now in order to synthesize a time-optimal control, it remains to compute the first switching time and switching point, since the consequent switching sequence can be constructed thereafter based on the procedure described above. Given an initial state $\Theta_0=(0,0)'$, the first switching time and point $p_1$ will be determined according to the target state, e.g., $\Theta_d=(2m_1\pi,2m_2\pi)'$, where $m_1,m_2\in\mathbb{Z}^+$, in such a way that the optimal trajectory follows a bang-bang control derived based on $p_1$ will reach $\Theta_d$. Under this construction, we may end up with a finite number of feasible trajectories starting with either $X$- or $Y$-trajectory, which reach the desired terminal state. The minimum time trajectory is then selected among them. 

\begin{figure}[t]
\subfigure[]{\includegraphics[scale=0.19]{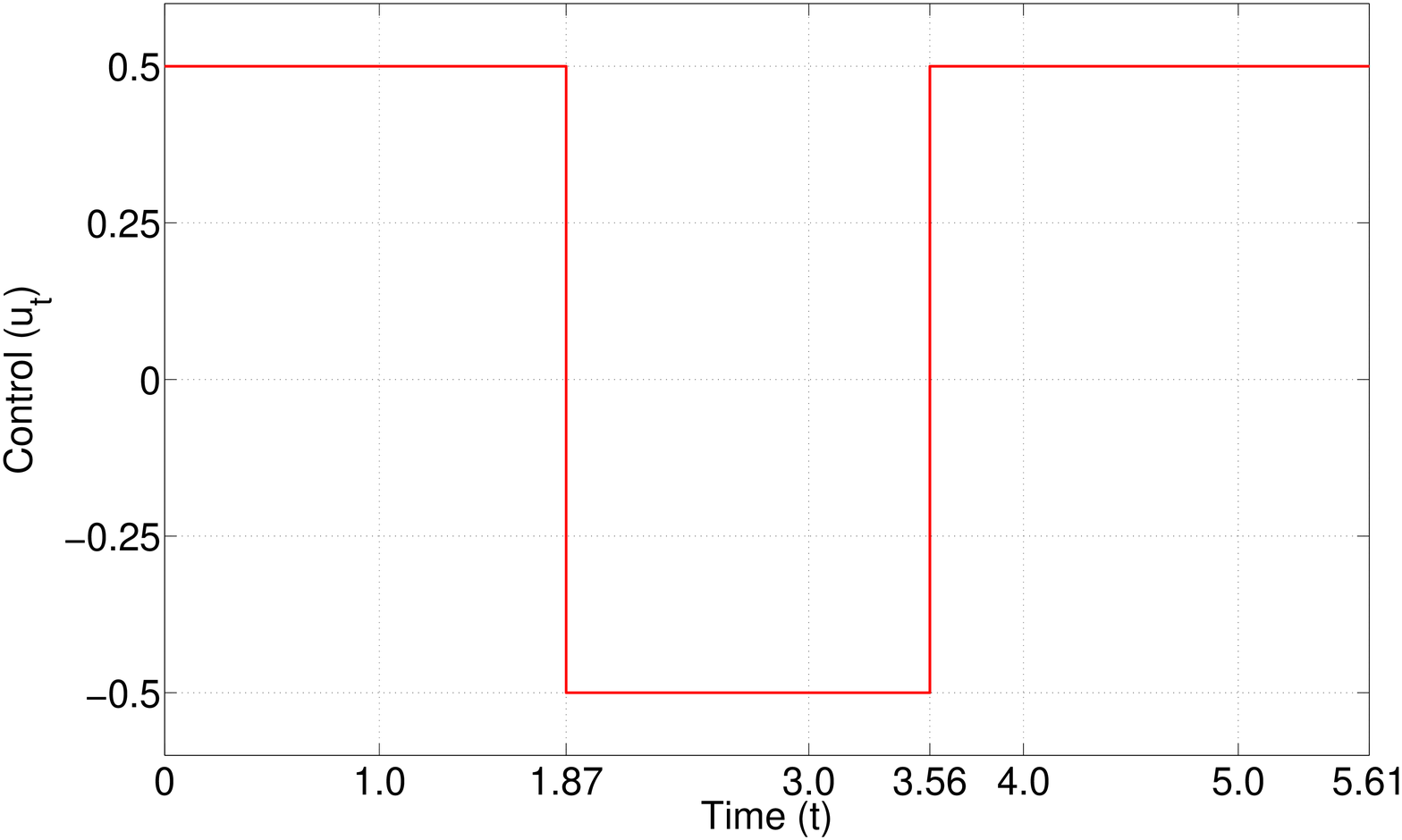}\label{fig:time_opt_con_2theta_neuron}}
\subfigure[]{\includegraphics[scale=0.19]{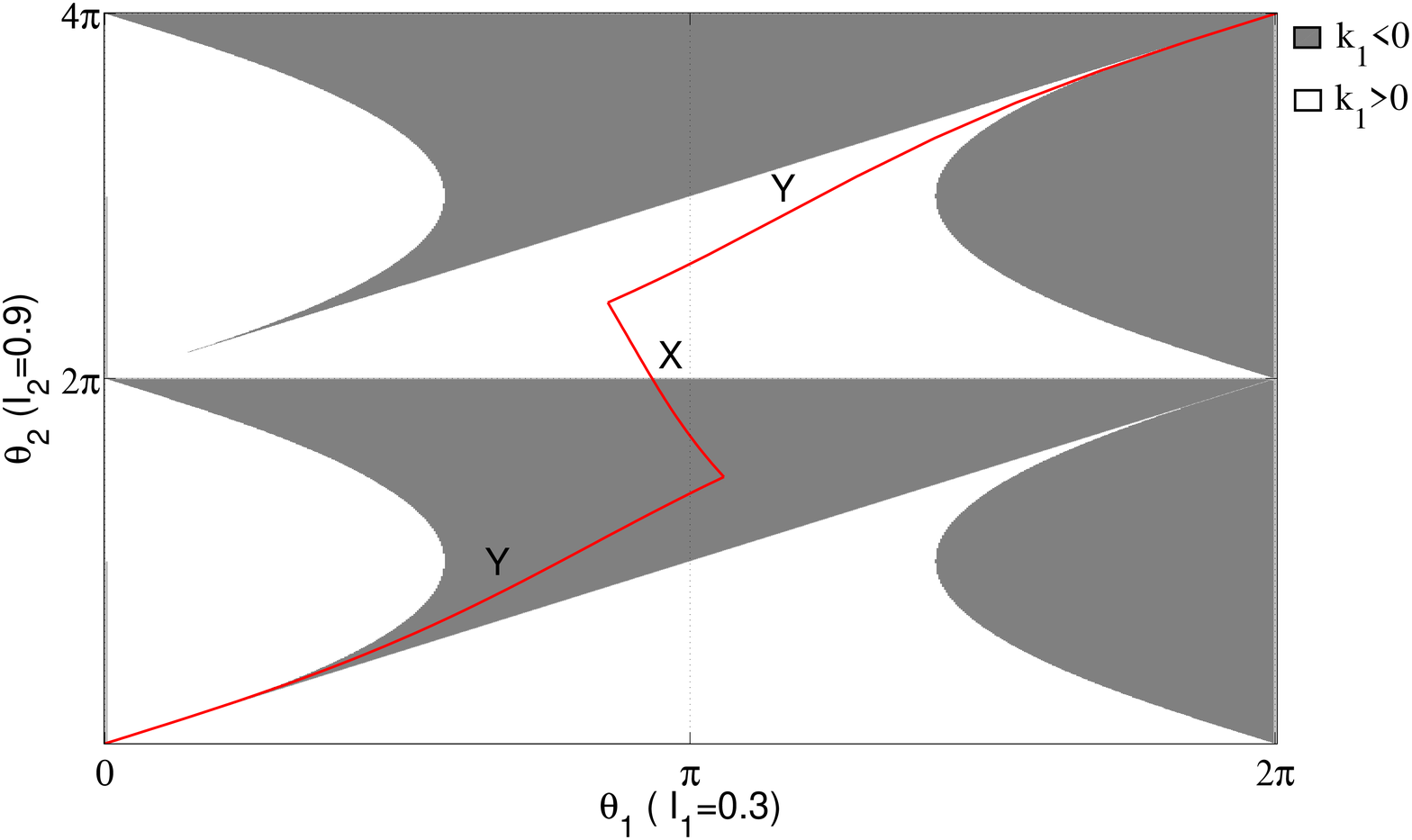}\label{fig:time_opt_traj_2theta_neuron}}
 \caption{\subref{fig:time_opt_con_2theta_neuron}  Time optimal control for two Theta neuron system with $I_1=0.3\ (\alpha_1=1.3,\,\beta_1=0.7)$ and $I_2=0.9\ (\alpha_2=1.9,\,\beta_2=0.1)$ to reach ($2\pi,4\pi$) with the control bounded by $M=0.5$ and  \subref{fig:time_opt_traj_2theta_neuron} corresponding trajectories. The gray and white regions represent where $k_1$ is negative and positive, respectively.}
	\label{fig:double_theta_neuron}
\end{figure}

Fig. \ref{fig:double_theta_neuron} illustrates an example of driving two Theta neurons time-optimally from $(0,0)'$ to $(2\pi,4\pi)'$ with the control bound $M=0.5$, where the natural frequencies of the oscillators are $\omega_1=1.1\ (I_1=0.3)$ and $\omega_2=1.9\ (I_2=0.9)$ corresponding to  $\alpha_1=1.3$, $\beta_1=0.7$ and $\alpha_2=1.9$, $\beta_2=0.1$. In this example, the time-optimal control has two switches at $t=1.87$ and $t=3.56$ and the minimum time is $5.61$.

% =============== Figure 3 ==================

% =================================================
\subsection{Simultaneous Control of Neuron Ensembles} \label{sec:simultaneous}
The complexity of deriving optimal controls for higher dimensional systems, i.e., more than two neurons, grows rapidly, and it makes sense to find out how the control of two neurons relates to the control of many. One may wonder whether it is possible to use a (optimal) control that spikes two neurons to manipulate an ensemble of neurons whose natural frequencies lie between those of the two nominal systems. Of course, if trajectories of the neurons with different frequencies have no crossings following a common control input, then the control designed for any two neurons guarantees to bound trajectories of all the neurons with their frequencies within the range of these two nominal neurons, whose trajectories can then be thought of as the envelope of these other neuron trajectories. We now show that this is indeed the case. 

\begin{lem}
	\label{lem:crossing}
	The trajectories of any two Theta neurons with positive baseline currents following a common control input have no crossing points.

	{\it Proof.} Consider two Theta neurons modeled by
	\begin{align}	
		\label{eq:theta1}
		\dot{\t_1} &= (1+I_1)+(1-I_1)\cos\t_1+(1-\cos\t_1)u,\quad \t_1(0)=0,\\
		\label{eq:theta2}
		\dot{\t_2} &= (1+I_2)+(1-I_2)\cos\t_2+(1-\cos\t_2)u, \quad \t_2(0)=0,
	\end{align}
	with positive baseline currents, $I_1,I_2>0$, and assume that $\w_1<\w_2$, which implies $I_1<I_2$ since $I_i=\frac{\w_i^2}{4}$, $i=1,2$. In the absence of any control input, namely, $u=0$, it is obvious that $\t_1(t)<\t_2(t)$ for all $t>0$ since $I_1<I_2$. Suppose that $\t_1(t)<\t_2(t)$ for $t\in(0,\tau)$ and these two phase trajectories meet at time $\tau$, i.e., $\t_1(\tau)=\t_2(\tau)$. Then, we have
$\dot{\t_1}(\tau)-\dot{\t_2}(\tau)=(I_1-I_2)(1-\cos(\t_1(\tau)))\leq 0$ and the equality holds only when the neurons spike at time $\tau$, i.e., $\t_1(\tau)=\t_2(\tau)=2n\pi$, $n\in\mathbb{Z}^+$. As a result, $\t_1(\tau^+)<\t_2(\tau^+)$, because $\t_1(\tau)=\t_2(\tau)$ and $\dot{\t_1}(\tau)<\dot{\t_2}(\tau)$, and hence there exist no crossings between the two trajectories $\t_1(t)$ and $\t_2(t)$. \hfill$\Box$
\end{lem}

Note that the same result as Lemma \ref{lem:crossing} holds and can be shown in the same fashion for both Sinusoidal and SNIPER phase models, as described in Section \ref{sec:sniper} and \ref{sec:sinusoidal}, when the model-dependent constant $z_1>z_2$ if $\w_1<\w_2$, which is in general the case. For example, in the SNIPER phase model, $z$ conventionally takes the form $z=2/\w$ as presented in Section \ref{sec:sniper}.

This critical observation extremely simplifies the design of external stimuli for spiking a neuron ensemble with different oscillation frequencies based on the design for two neurons with the extremal frequencies over this ensemble. We illustrate this important result by designing optimal controls for two Theta and two Sinusoidal neurons employing the Legendre pseudospectral method, which will be presented in Section \ref{sec:ps}. Fig. \ref{fig:twoneuron} shows the optimized controls and corresponding trajectories for Theta and Sinusoidal neurons with their frequencies $\w$ belonging to $[0.9,1.1]$ and $[1.0,1.1]$, respectively. The optimal controls are designed based only on the extremal frequencies of these two ranges, i.e., 0.9 and 1.1 for the Theta neuron model and 1.0 and 1.1 for the Sinusoidal model.
%, and for the cost functional
%$$J = \|\Theta_d-\Theta(T)\|^2 + \alpha \int_0^{T} u^2(t)dt\ d\w,$$
%which minimizes the terminal error and input energy with a relative scaling given by $\alpha$, where $\Theta_d=(2\pi,2\pi)'$, and we choose $\alpha = 0.1$ and $T=2\pi$. The controlled (black) and uncontrolled (gray) state trajectories clearly show the improvement in simultaneous spiking of the ensemble of neurons. While a bound is necessary to provide in practice, the inclusion of the minimum energy term in the cost function serves to regularize the control against high amplitude values.

This design principle greatly reduces the complexity of finding controls to spike a large number of neurons. Although the optimal control for two neurons is in general not optimal for the others, this method produces a good approximate optimal control. In the next section, we will introduce a multivariate pseudospectral computational method for constructing optimal spiking or synchronization controls.

% This critical observation extremely simplifies the design of external stimuli for spiking a neuron ensemble with different oscillation frequencies based on the design for two neurons with the extremal frequencies over this ensemble. We illustrate this important result by applying the above time-optimal control for two neurons to the neurons with their frequencies within the range of these two. The phase trajectories following this time-optimal control are shown in Figure xxx. This design principle greatly reduces the complexity of finding controls to spike a large number of neurons. Although the optimal control for two neurons is in general not optimal for the others, this method produces a good approximate optimal control. In the next section, we will introduce a multivariate pseudospectral computational method for constructing optimal spiking controls.

%%%%%%%%%%%%%%%%%%%%%%% Pseudospectral Method (Many Neurons) %%%%%%%%%%%%%%%%%%%%%%%%
\section{Computational Optimal Control of Spiking Neuron Networks}
\label{sec:ps}
As we move to consider the synthesis of controls for neuron ensembles, the analytic methods used in the one and two neuron case become impractical to use. As a result, developing computational methods to derive inputs for ensembles of neurons is of particular practical interest. We solve the optimal control problem in \eqref{eq:steering} using a modified pseudospectral method. Global polynomials provide accurate approximations in such a method which has shown to be effective in the optimal ensemble control of quantum mechanical systems \cite{Li_PNAS11, Ruths_JCP11, Li_IEEE_QCP11, Li_CDC11_PS}. Below we outline the main concepts of the pseudospectral method for optimal control problems and then show how it can be extended to consider the ensemble case.

% As we move to consider the synthesis of controls for neuron ensembles, the analytic methods used in the one and two neuron case become impractical to use. As a result, developing computational methods to derive inputs for ensembles of neurons is of particular practical interest. We solve the optimal control problem in \eqref{eq:steering} and the optimal ensemble control problem in \eqref{eq:neuronensemble} using a modified pseudospectral method. Global polynomials provide accurate approximations in such a method, which is motivated by the critical role played by polynomials in the proof of ensemble controllability (see Section \ref{sec:neuron}). Such a method has shown to be effective in the optimal ensemble control of quantum mechanical systems \cite{Li_PNAS11, Ruths_JCP11}. Below we outline the main concepts of the pseudospectral method for optimal control problems and then show how it can be extended to consider the ensemble case.

Spectral methods involve the expansion of functions in terms of orthogonal polynomial basis functions on the domain $[-1,1]$ (similar to Fourier series expansion), facilitating high accuracy with relatively few terms \cite{Canuto06}. The pseudospectral method is a spectral collocation method in which the differential equation describing the state dynamics is enforced at specific nodes. Developed to solve partial differential equations, these methods have been recently adopted to solve optimal control problems \cite{Elnagar95, Ross03, Fahroo01}. We focus on Legendre pseudospectral methods and consider the transformed optimal control problem on the time domain $[-1,1]$.

The fundamental idea of the Legendre pseudospectral method is to approximate the continuous state and control functions, $\Theta(t)$ and $u(t)$, by $N^{th}$ order Lagrange interpolating polynomials, $I_N \Theta(t)$ and $I_N u(t)$, based on the Legendre-Gauss-Lobatto (LGL) quadrature nodes, which are defined by the union of the endpoints, $\{-1,1\}$, and the roots of the derivative of the $N^{th}$ order Legendre polynomial. Note that the non-uniformity in the distribution of the LGL nodes and the high density of nodes near the end points are a key characteristic of pseudospectral discretizations by which the Runge phenomenon is effectively suppressed \cite{Fornberg98}. The interpolating approximations of the state and control functions, $\Theta(t)$ and $u(t)$ can be expressed as functions of the Lagrange polynomials, $\ell_k(t)$, \cite{Szego59}
\begin{align}
	\label{eq:Ix}
	\Theta(t) &\approx I_N \Theta(t) = \sum_{k=0}^N \bar{\Theta}_k \ell_k(t),\\
u(t) &\approx I_N u(t) = \sum_{k=0}^N \bar{u}_k \ell_k(t).\nonumber
\end{align}
The derivative of $I_N \Theta(t)$ at the LGL node $t_j$, $j=0,1,\ldots,N$, is given by
\begin{equation}\label{eq:dinterpsc_j}
	\frac{d}{dt} I_N \Theta(t_j) = \sum_{k=0}^N \bar{\Theta}_k \dot{\ell}_k(t_j)=\sum_{k=0}^N D_{jk}\bar{\Theta}_k, \nonumber
\end{equation}
where $D_{jk}$ are elements of the constant $(N+1)\times(N+1)$ differentiation matrix \cite{Canuto06}. Finally, the integral cost functional in the optimal control problem \eqref{eq:steering} can be approximated by the Gauss-Lobatto integration rule, and we ultimately convert the optimal control problem into the following finite-dimensional constrained minimization problem
\begin{align}
	\min\ \ & \frac{T}{2}\sum_{j=0}^N \bar{u}_j^2\, w_j \nonumber\\
	\label{eq:dps}{\rm s.t.}\ \ & \sum_{k=0}^N D_{jk} \bar{\Theta}_k =\frac{T}{2}\Big[f(\bar{\Theta}_j)+\bar{u}_j\, Z(\bar{\Theta}_j)\Big],\\
	& \bar{\Theta}(-1)=0 ,\nonumber\\
	& \bar{\Theta}(1)=\Theta_d, \quad \forall\ j\in\{0,1,\ldots,N\},\nonumber
\end{align}
where $\Theta_d=(2m_1\pi,2m_2\pi,\ldots,2m_n\pi)'$, $m_i\in\mathbb{Z}^+$, $i=1,\ldots, n$, is the target state and $w_j$ are the LGL weights given by $$w_j=\frac{2}{N(N+1)}\frac{1}{(L_N(t_j))^2},$$
in which $L_N$ is the $N^{th}$ order Legendre polynomial. Solvers for this type of constrained nonlinear programs are readily available and straightforward to implement.

% ============= Remark 2 =================
\begin{rmk}[Extension to an infinite ensemble of neuron systems]
\label{rmk:continuum}
The pseudospectral computational method can be readily extended to consider an infinite population of neurons, for instance, with the frequency distribution over a closed interval, $\w\in[\w_a,\w_b]\subset\mathbb{R}^+$. In such a case, the parameterized state function can be approximated by a two-dimensional interpolating polynomial, namely, $\Theta(t,\w)\approx I_{N\times N_{\w}} \Theta(t,\w)$, based on the LGL nodes in the time $t$ and the frequency $\w$ domain. Similarly, the dynamics of the state can be expressed as an algebraic constraint and a corresponding minimization problem can be formed \cite{Ruths_JCP11}.

\end{rmk}

\subsection{Optimized Ensemble Controls}
\label{sec:oec}

% =================== Figure 4 ==================
\begin{figure}[t]
\centering
\begin{tabular}{cc}
	\includegraphics[width=0.5\linewidth]{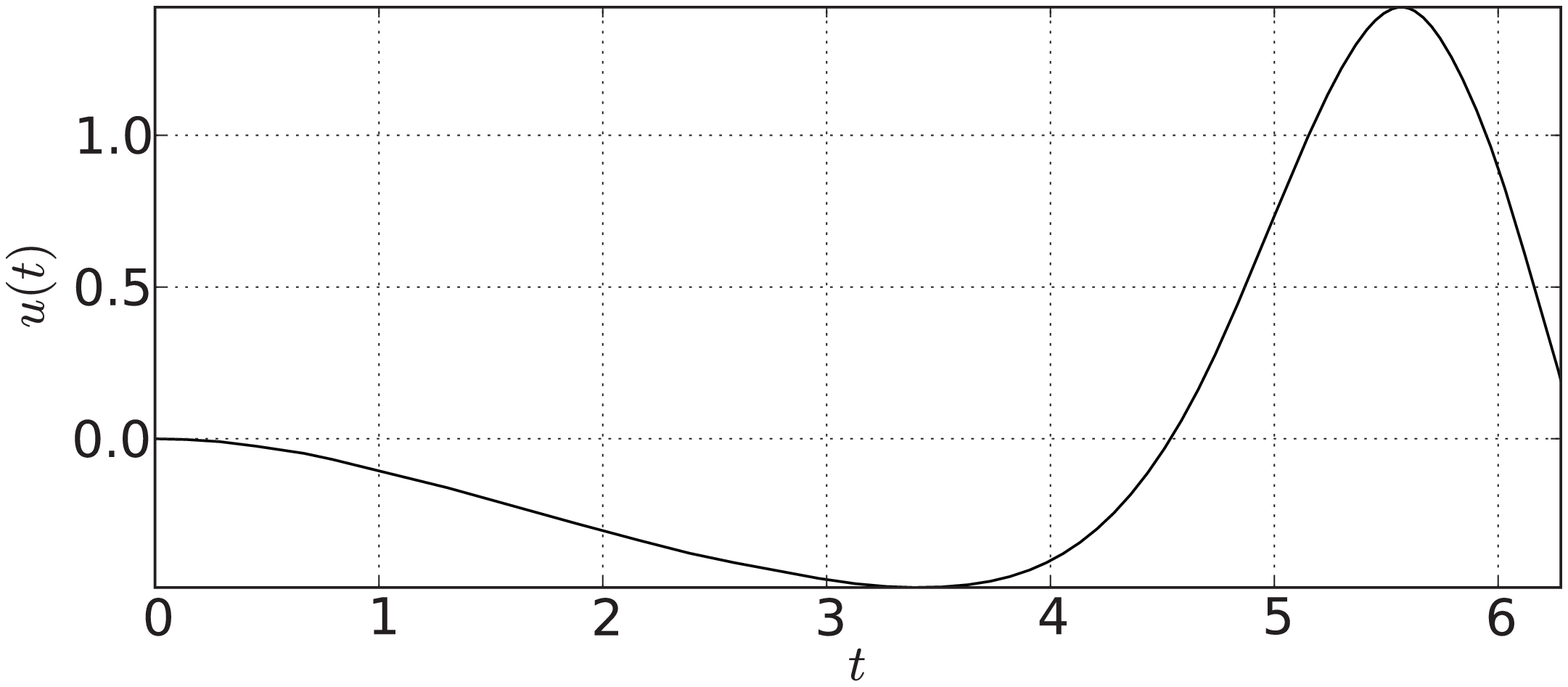} &
	\includegraphics[width=0.5\linewidth]{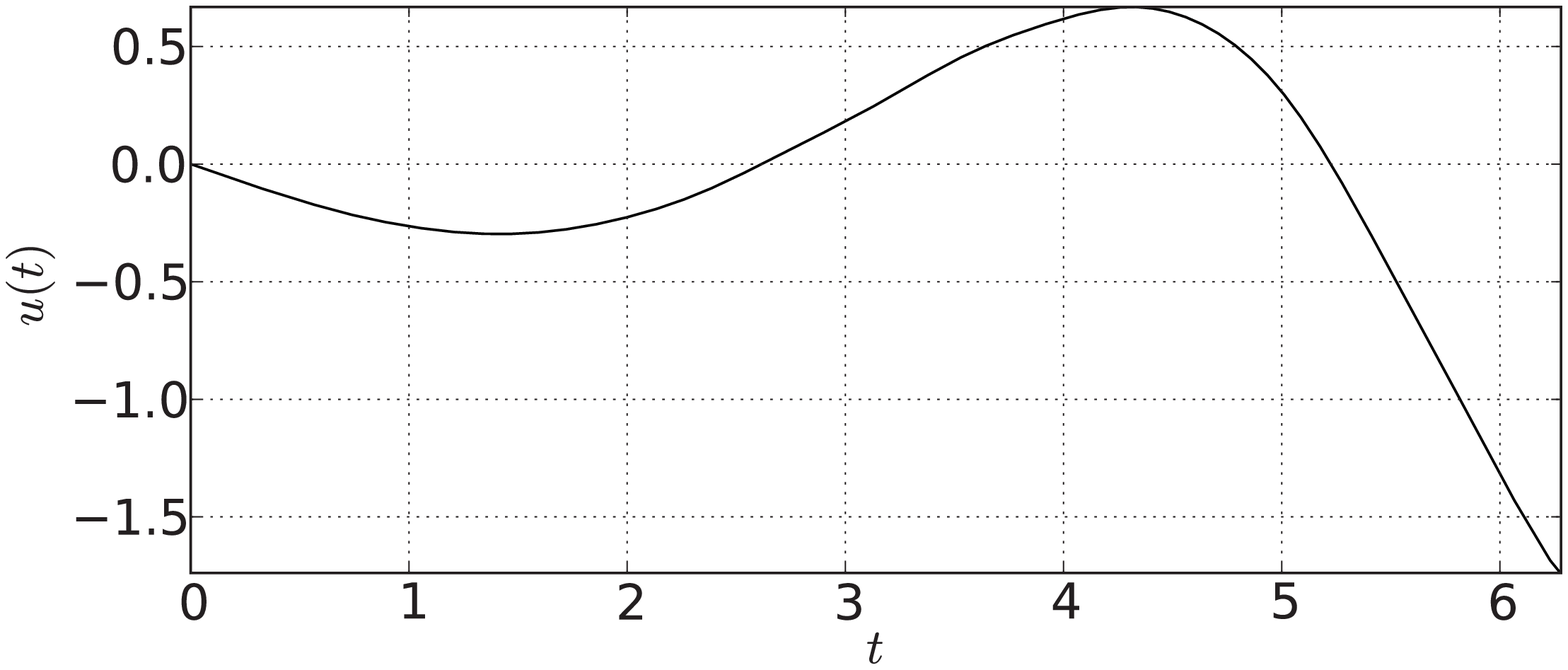} \\
	\includegraphics[width=0.5\linewidth]{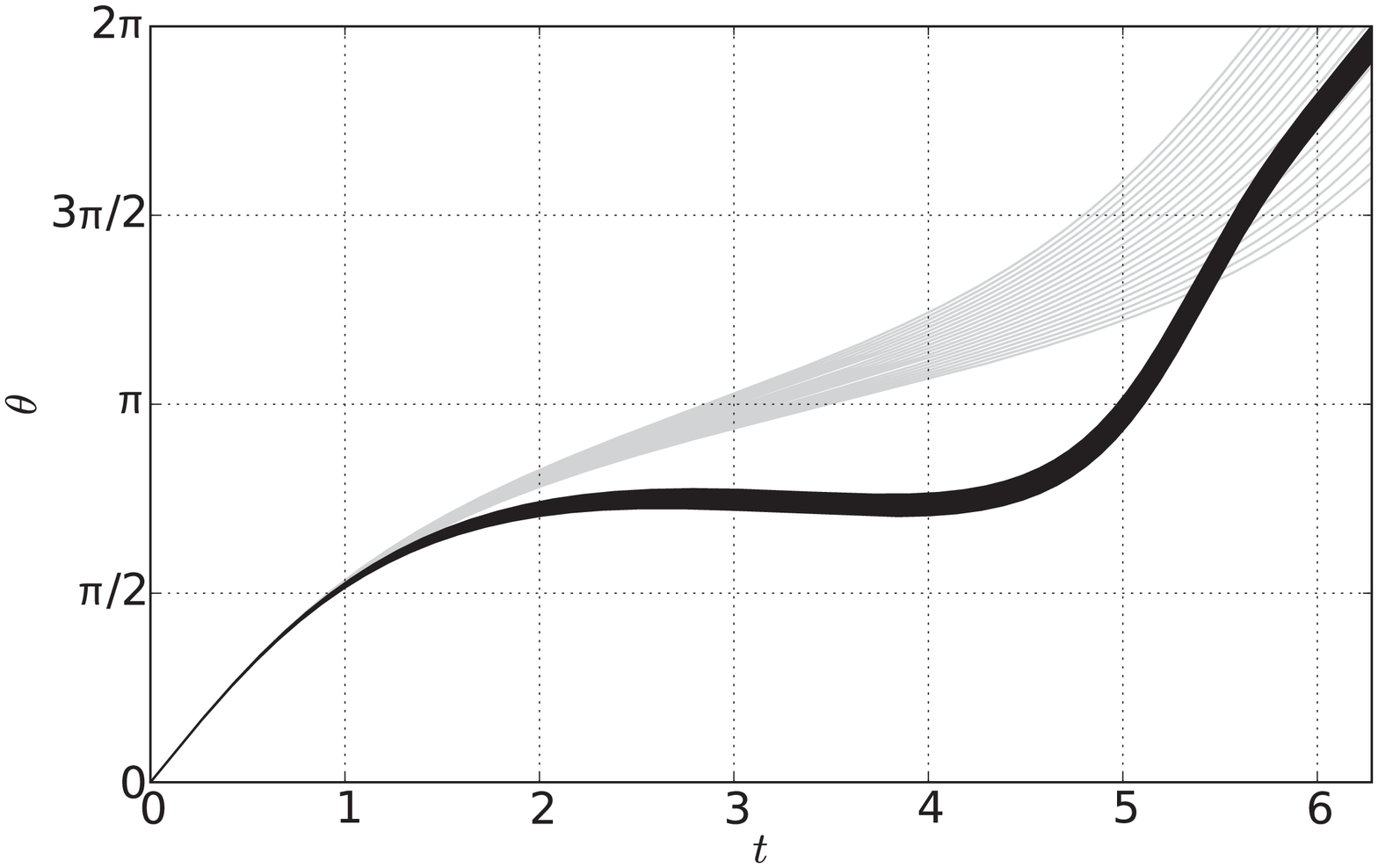} &
	\includegraphics[width=0.5\linewidth]{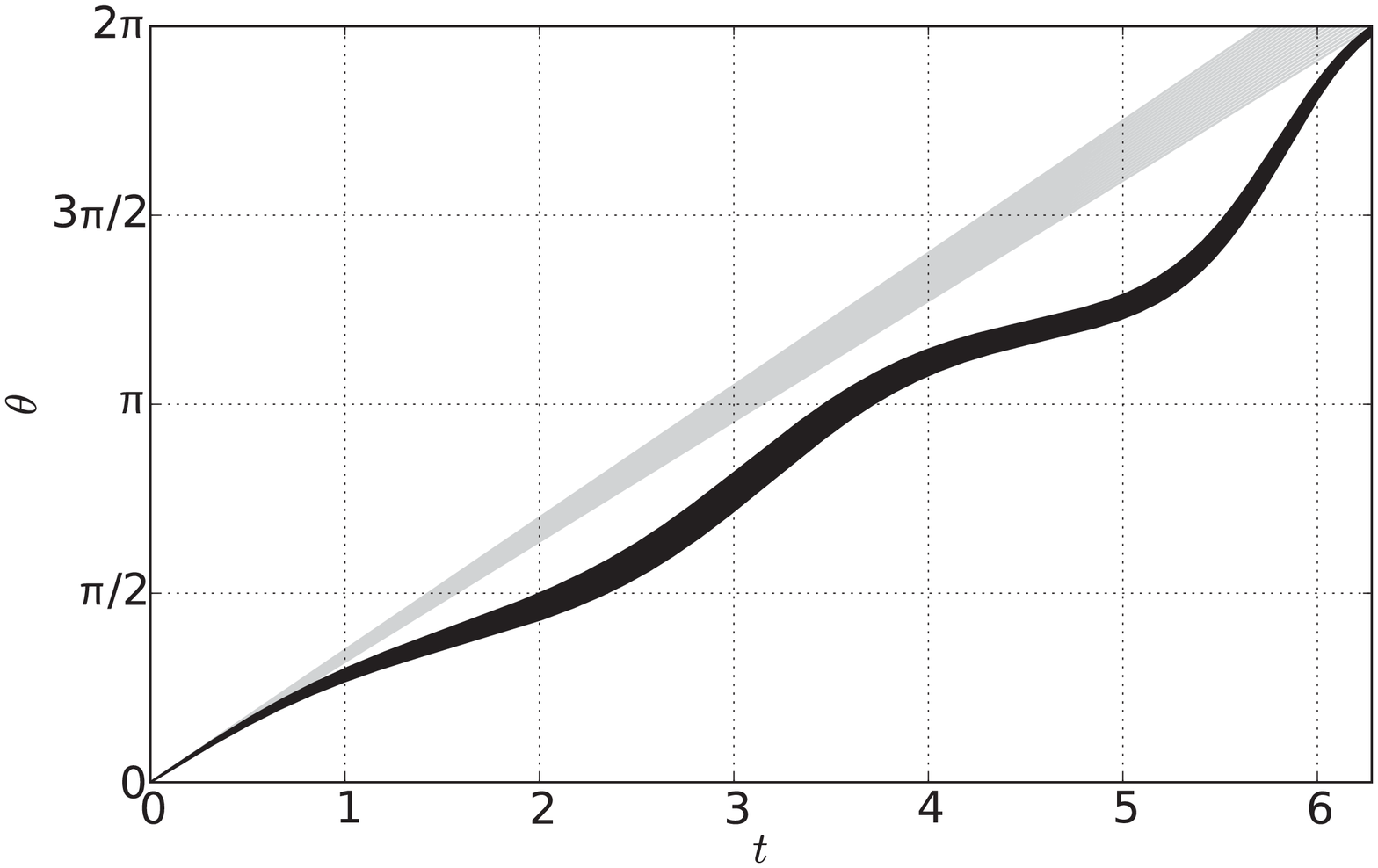}
\end{tabular}
\caption{The controls (top) and state trajectories (bottom) of Theta (left) and Sinusoidal (right) PRC neuron models (for $\alpha=1$, $\beta=0.1$, $T=2\pi$ in \eqref{eq:J}). The Theta model is optimized for $\w\in[0.9,1.1]$. The Sinusoidal model is optimized for $\w\in[1.0,1.1]$. The gray states correspond to uncontrolled state trajectories, and provide a comparison for the synchrony improvement provided by the compensating optimized ensemble control.} \label{fig:twoneuron}
\end{figure}

We can now apply the above methodology to synthesize optimal controls for neuron ensembles. Since neurons modeled by the SNIPER PRC are special cases of the Theta neuron, here we consider Theta and Sinusoidal neuron models. The computational method outlined above permits a flexible framework to optimize based on a very general cost functional subject to general constraints. We illustrate this by selecting an objective of the type,
\begin{eqnarray}
	\label{eq:J}
	J=\alpha\|\Theta_d-\Theta(T)\|^2 + \beta \int_0^{T} u^2(t)dt,
\end{eqnarray}
% \begin{equation}
% J = \alpha\int_{\w_a}^{\w_b} \big(2\pi-\Theta(t_1,\w)\big)^2 + \beta \int_0^{t_1} u^2(t)dt\ d\w,
% \end{equation}
which minimizes the terminal error and input energy with a relative scaling given by the constants $\alpha$ and $\beta$. In highly complex problems, such as those given by ensemble systems as described in Remark \ref{rmk:continuum}, this scaling provides a tunable parameter that determines the trade-off between performance and input energy. 

% =================== Figure 5 ==================
\begin{figure}[t]
\centering
\begin{tabular}{cc}
	\includegraphics[width=0.5\linewidth]{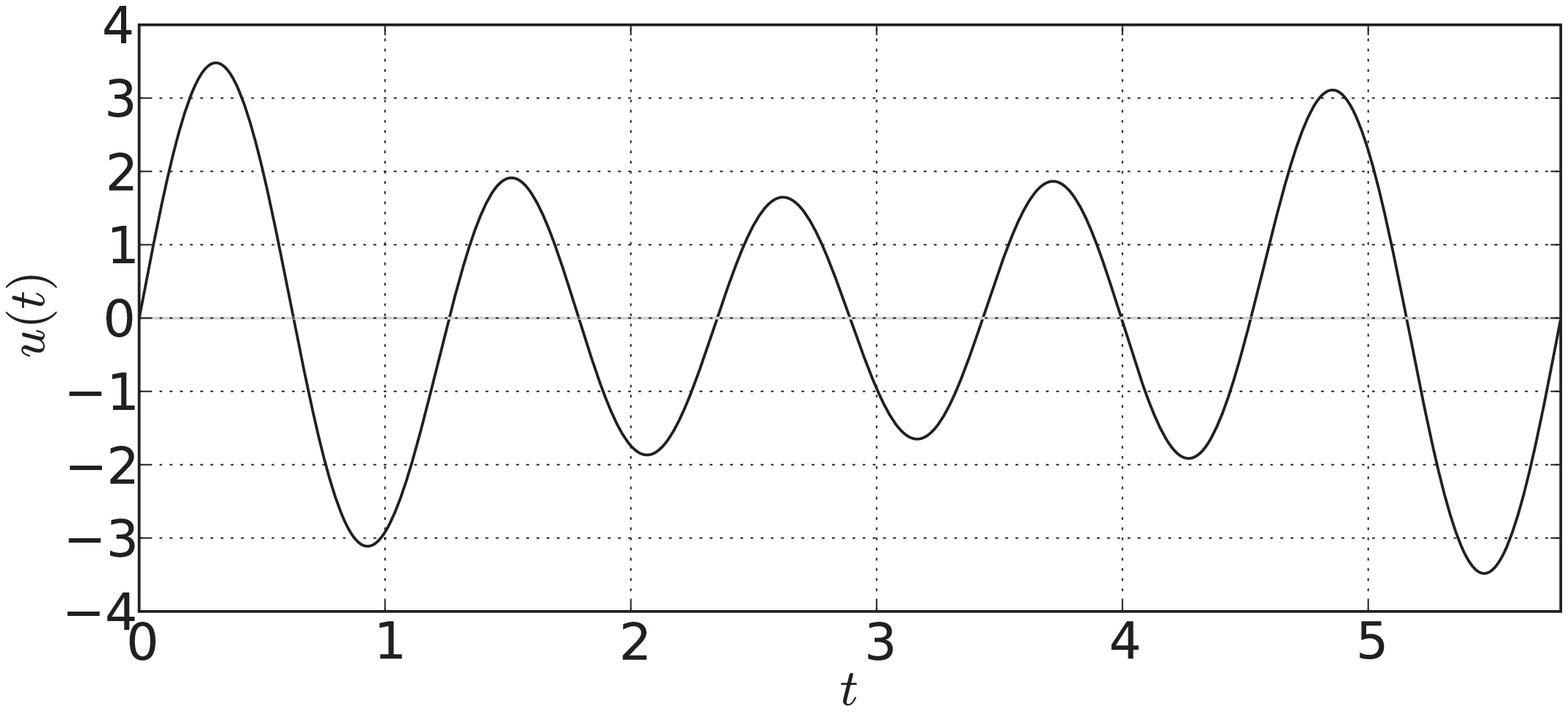} &
	\includegraphics[width=0.5\linewidth]{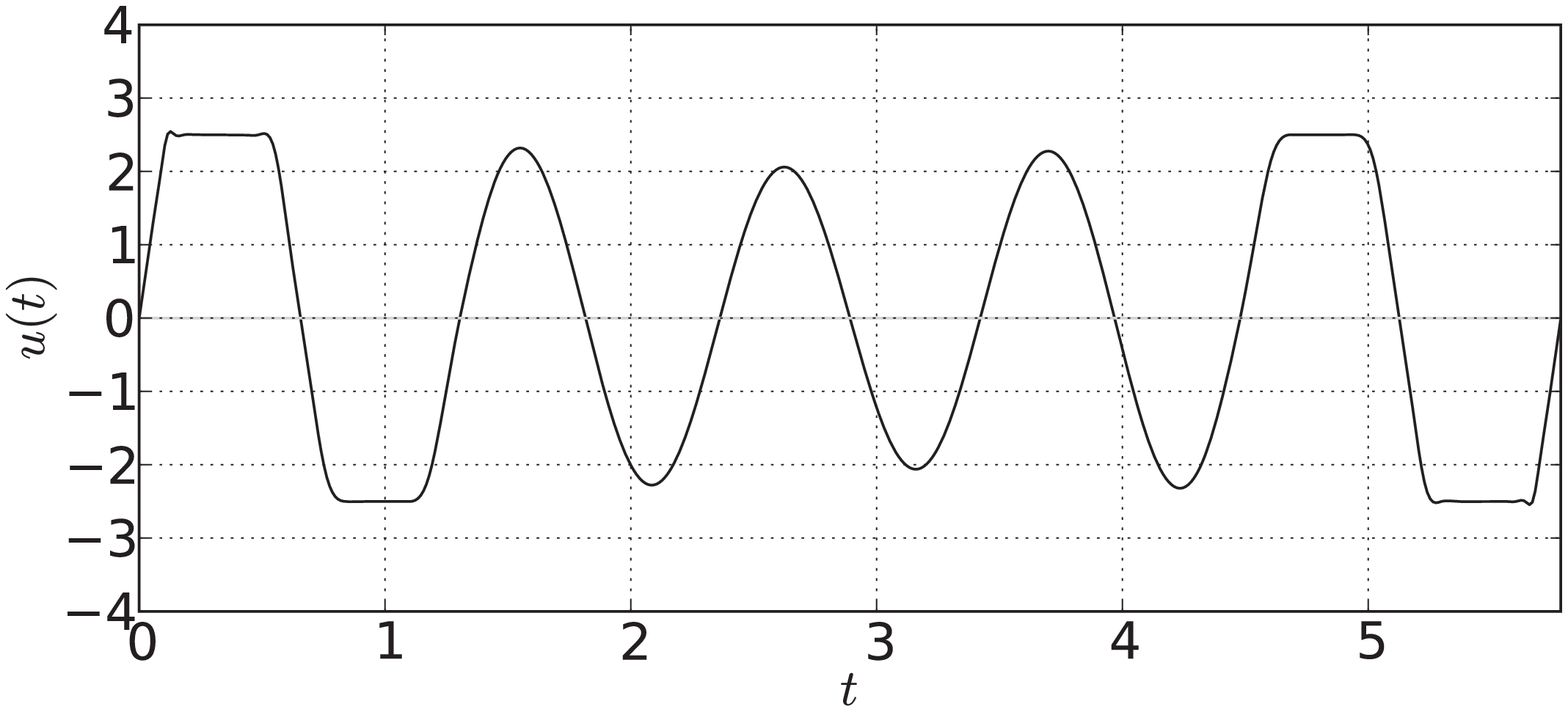} \\
	\includegraphics[width=0.5\linewidth]{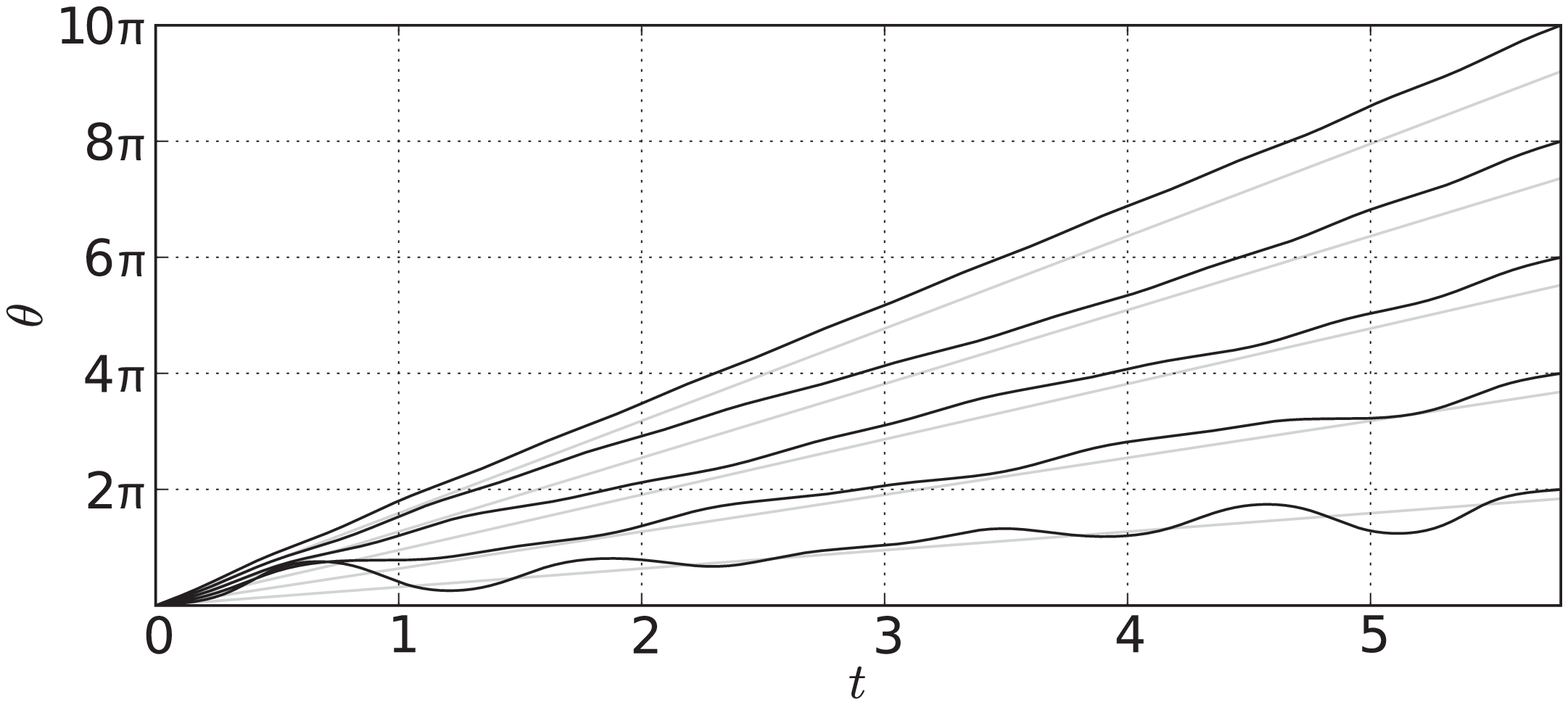} &
	\includegraphics[width=0.5\linewidth]{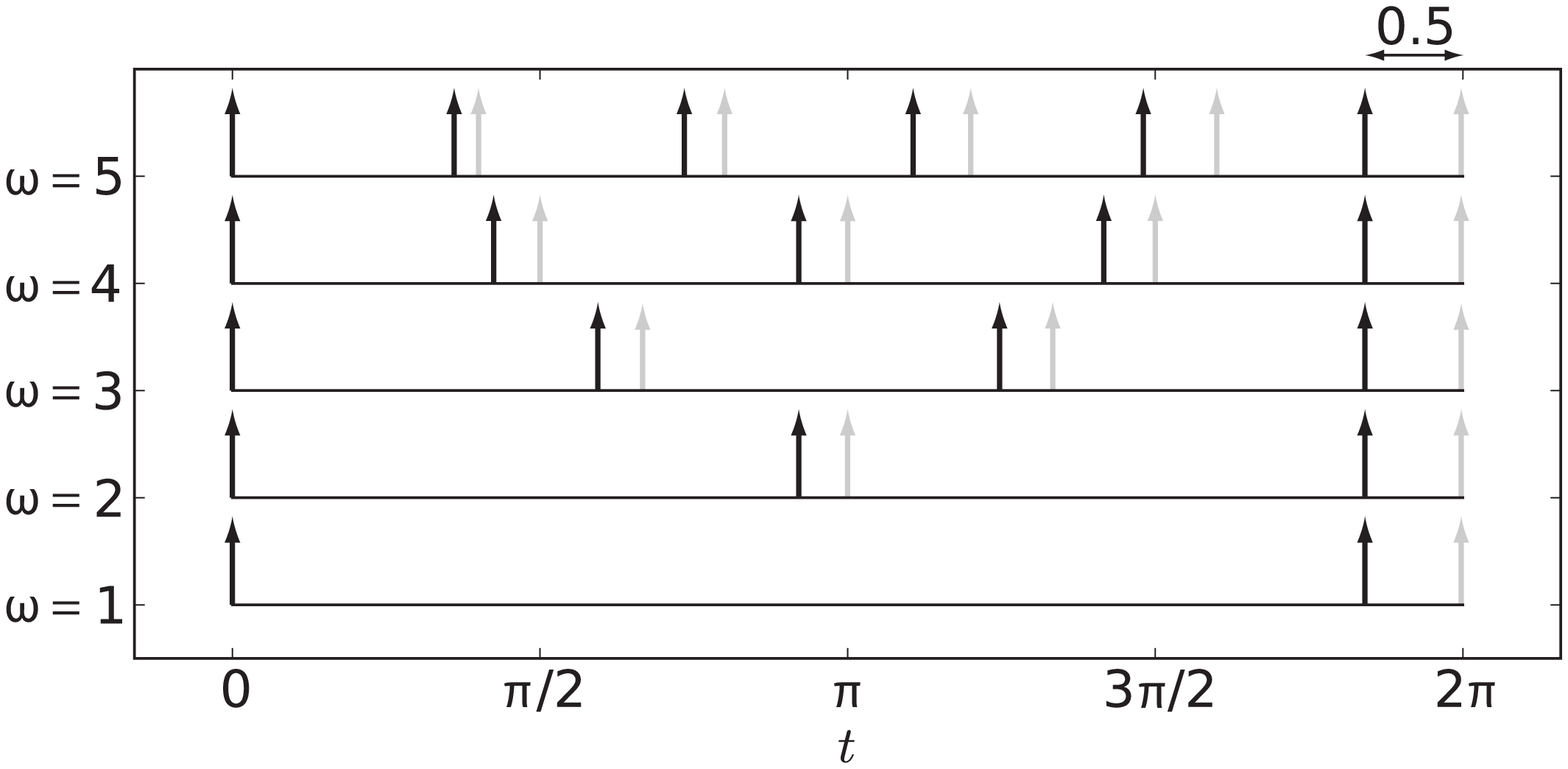}
\end{tabular}
\caption{The controls and amplitude constrained controls, $A\leq2.5$, upper left and right respectively, of a Sinusoidal PRC neuron model driving five frequencies, $(\w_1,\w_2,\w_3,\w_4,\w_5)=(1,2,3,4,5)$, to the desired targets $\Theta(T)=(2\pi,4\pi,6\pi,8\pi,10\pi)$ when $T=2\pi-0.5$. These controls yield highly similar state trajectories (left, shown for unconstrained control) and spiking sequences (right, shown for constrained control), which corresponds to when the state trajectories cross multiples of $2\pi$. Black coloring indicates a controlled state trajectory or spike sequence, whereas gray coloring indicates a trajectory or spike sequence without control.} \label{fig:sinusoidal}
\end{figure}

% =================== Figure 6 ==================
\begin{figure}[t]
\centering
\begin{tabular}{cc}
	\includegraphics[width=0.5\linewidth]{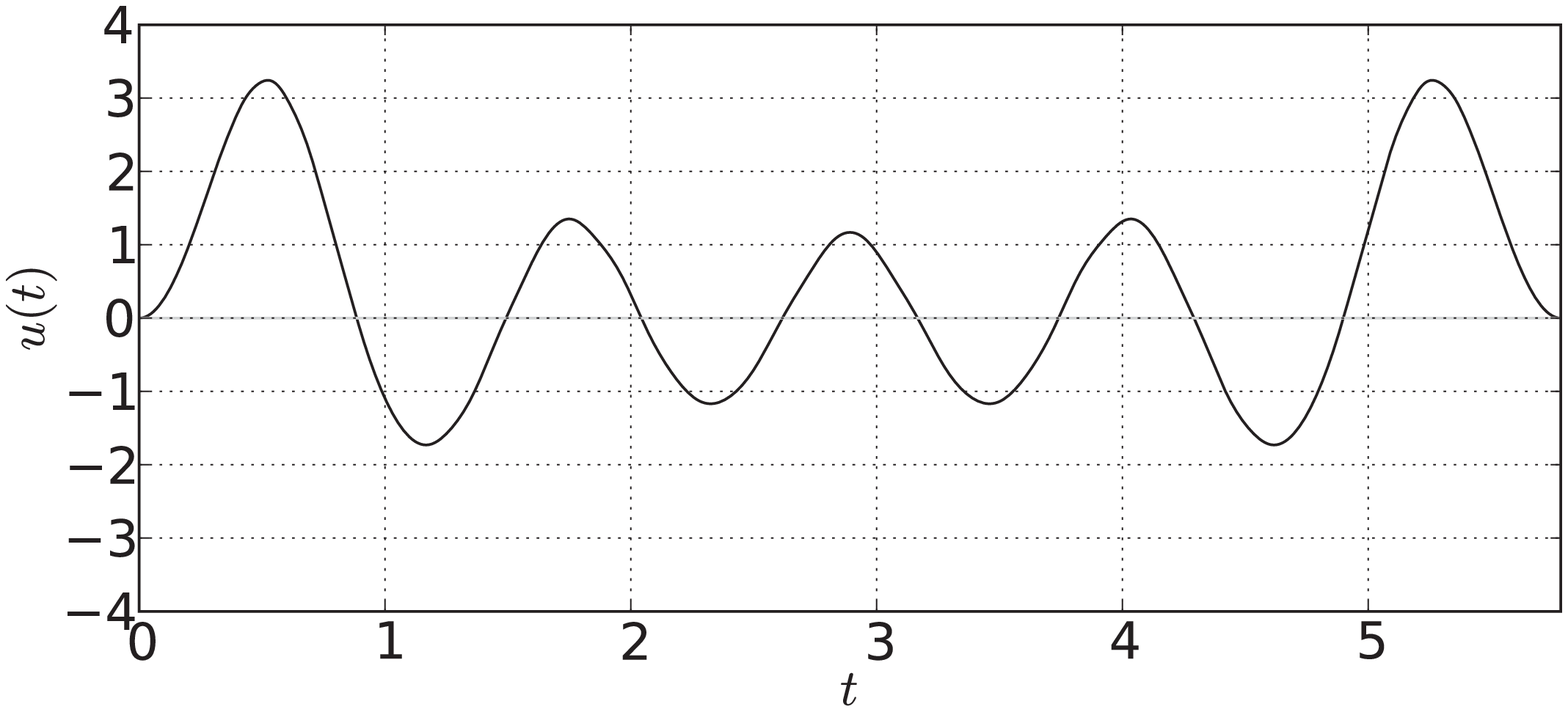} &
	\includegraphics[width=0.5\linewidth]{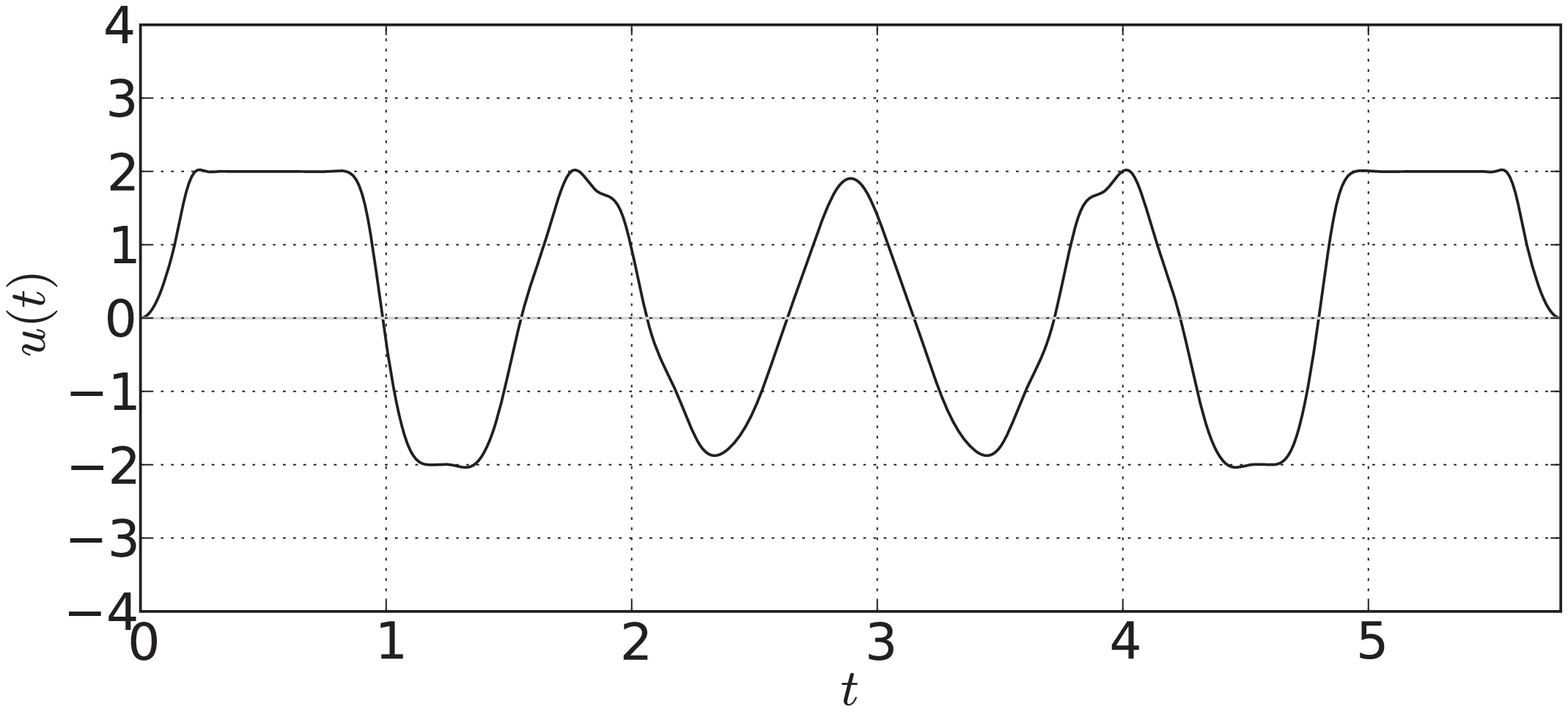} \\
	\includegraphics[width=0.5\linewidth]{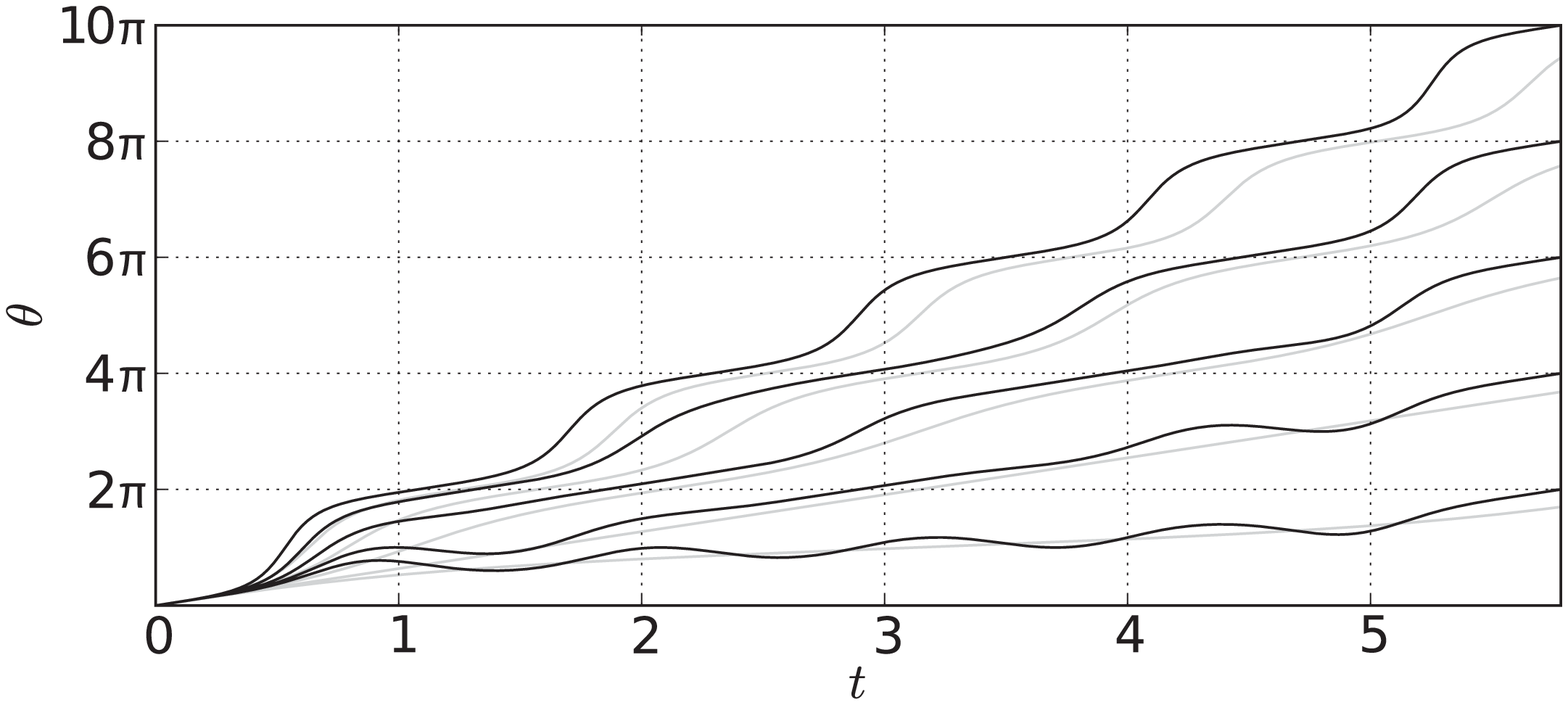} &
	\includegraphics[width=0.5\linewidth]{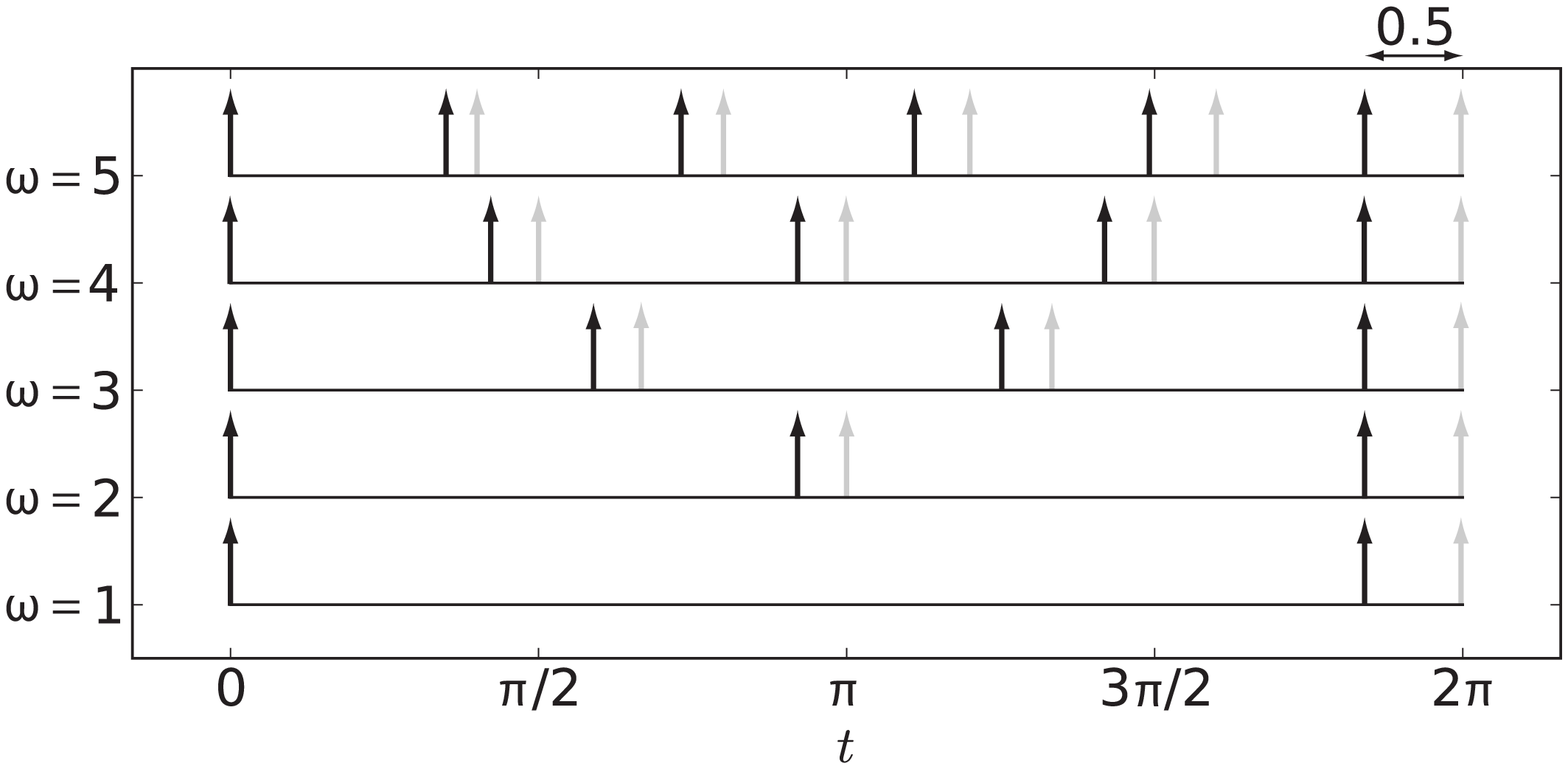}
\end{tabular}
\caption{The controls and amplitude constrained controls, $A\leq2$, upper left and right respectively, of a Theta PRC neuron model driving five frequencies, $\mathbf{\w}=(\w_1,\w_2,\w_3,\w_4,\w_5)=(1,2,3,4,5)$, to the desired targets $\Theta(T)=(2\pi,4\pi,6\pi,8\pi,10\pi)$ when $T=2\pi-0.5$. These controls yield highly similar state trajectories (left, shown for unconstrained control) and spiking sequences (right, shown for constrained control), which corresponds to when the state trajectories cross multiples of $2\pi$. Black coloring indicates a controlled state trajectory or spike sequence, whereas gray coloring indicates a trajectory or spike sequence without control.} \label{fig:theta}
\end{figure}

Fig. \ref{fig:twoneuron} shows the optimized controls and corresponding trajectories for Theta and Sinusoidal neuron models for $\alpha=1$, $\beta = 0.1$, $T=2\pi$, and $\w$ belongs to $[0.9,1.1]$ and $[1.0,1.1]$ respectively. In this optimization, the controls are optimized over the two neuron systems with extremal frequencies, whose trajectories form an envelope, bounding the trajectories of other frequencies in between, as described in Section \ref{sec:simultaneous}. We are able to design compensating controls for the entire frequency band solely by considering these upper and lower bounding frequencies. The controlled (black) and uncontrolled (gray) state trajectories clearly show the improvement in simultaneous spiking of the ensemble of neurons. While a bound is necessary to provide in practice, the inclusion of the minimum energy term in the cost function serves to regularize the control against high amplitude values.

In Fig. \ref{fig:sinusoidal}, we demonstrate the flexibility of the method to drive multiple Sinusoidal neurons to desired targets. In particular we seek to simultaneously spike five frequencies with widely dispersed frequency values at a time $T$ different from their natural period. In this figure we consider the frequencies $(\w_1,\w_2,\w_3,\w_4,\w_5)=(1,2,3,4,5)$ and design controls to drive these systems to $(2\pi,4\pi,6\pi,8\pi,10\pi)$, respectively, at a time $T=2\pi-0.5$. Controls for minimum energy ($\alpha=0$, $\beta=1$) transfer can be designed for both the unconstrained and amplitude constrained cases (shown for $A\leq 2.5$). In both cases, the state trajectories and spike sequence (shown in the lower half of the figure) follow the same general pattern. The spike train shows that the controls are able to advance the firing of each neuron so that all five spike simultaneously at the desired terminal time. Again the gray coloring indicates uncontrolled trajectories or spike trains and offers a comparison of improvement in synchrony. 

Similarly, Fig. \ref{fig:theta} provides the same presentation as above for the minimum energy transfer for Theta neurons of the same frequencies to the same desired targets. In this case the constrained control is limited to $A\leq 2$.

%%%%%%%%%%%%%%%%%%%%%%%%%%%%%%%%%%%%%%%%%%%%%%%%%%%%%%%%%%%%%%%%%%%%%%%%%%%%%%%%%%%%%%%%%%%%
\section{Conclusion}
In this paper, we considered the control and synchronization of a neuron ensemble described by phase models. We showed that this ensemble system is controllable for various commonly-used phase models. We also derived minimum-power and time-optimal controls for single and two neuron systems. The development of such optimal controls is of practical importance, for example, in therapeutic procedures such as deep brain stimulation for Parkinson's disease and cardiac pacemakers for heart disease. In addition, we adopted a computational pseudospectral method for constructing optimal controls that spike neuron ensembles, which demonstrated the underlying controllability properties of such neuron systems. The methodology resulting from this work can be applied not only to neuron oscillators but also to any oscillating systems that can be represented using similar model reduction techniques such as biological, chemical, electrical, and mechanical oscillators. A compelling extension of this work is to consider networks of coupled oscillators, whose interactions are characterized by a coupling function, $H$, acting between each pair of oscillators. For example, in the well-known Kuramoto's model, the coupling between the $(i,j)$-pair is characterized by the sinusoidal function of the form $H(\t_i,\t_j)=\sin(\t_i-\t_j)$ \cite{Kuramoto84}. The procedure presented in Theorem \ref{thm:thetaneuron} can be immediately applied to examine controllability of interacting oscillators by investigating the recurrence properties of the vector field $f+H$, and the Lie algebra $\{f+H,Z\}_{LA}$. Similarly, the pseudospectral method presented in Section \ref{sec:ps} and its extension addressed in Remark \ref{rmk:continuum} can be employed to calculate optimal controls for spiking or synchronizing networks of coupled neurons with their natural frequencies varying on a continuum.

% The theoretical results presented in this paper characterize the fundamental limit of how the dynamics of neurons can be perturbed by the use of external inputs.

%\renewcommand{\appendixname}{Appendix \thesection}
%%%%%%%%%%%%%%%%%%%%%%%%%%%%%%%%%%%%%%%%%% Appendix %%%%%%%%%%%%%%%%%%%%%%%%%%%%%%%%%%%%%
% \appendix
\appendices

% ======================= Appendix A: Chow's Theorem =========================
\section{Chow's Theorem}    % Each appendix must have a short title.
\label{apd:Chow}
\begin{thm} (Versions of Chow's Theorem) Let $\{f_1(x),f_2(x),\ldots,f_m(x)\}$ be a collection of vector fields such that the collection
$\{f_1(x),f_2(x),\ldots,f_m(x)\}_{LA}$ is
\begin{enumerate}
    \item[a)] analytic on an analytic manifold $M$. Then given any point $x_0\in M$, there exists a maximal submanifold $N\subset M$ containing $x_0$ such that $\{\exp\{x_i\}\}_G\, x_0=\{\exp\{x_i\}_{LA}\}_{G}\, x_0=N$.
    \item[b)] $C^{\infty}$ on a $C^{\infty}$ manifold $M$ with dim (span$\{f_i(x)\}_{LA}$) constant on $M$. Then given any point $x_0\in M$, there exists a maximal submanifold $N\subset M$ containing $x_0$ such that $\{\exp\{x_i\}\}_G\, x_0=\{\exp\{x_i\}_{LA}\}_{G}\, x_0=N$.
\end{enumerate}
For more details, please see \cite{Brockett76}.
\end{thm}

% ======================== Appendix B: Spiking a Theta Neuron ====================
\section{Optimal Control of a Single Theta Neuron}

% =========== Theta neuron (unbounded control) ============
\subsubsection{Unbounded Minimum-Power Control of a Theta Neuron}
\label{apd:thetaneuron_unbounded}
The minimum-power control of a single Theta neuron is formulated as
\begin{align*}
	\min\quad & \int_0^{T}{u^2(t)}dt,\\
	{\rm s.t.}\quad & \dot{\theta}=\alpha+\beta\cos\theta+(1-\cos\theta)u(t),\\
	&\theta(0)=0, \quad\quad \theta(T)=2\pi.
\end{align*}
We then can form the control Hamiltonian,
\begin{equation}
    \label{eq:hamiltonian}
	H=u^2+\lambda(\alpha+\beta\cos\theta+u-u\cos\theta),
\end{equation}
where $\lambda$ is the Lagrange multiplier. The necessary conditions for optimality from the maximum principle yield
\begin{align}
	\label{eq:lambda_dot}
	\dot{\lambda} &= -\frac{\partial H}{\partial \theta}=\lambda(\beta-u)\sin\theta,\\
	\frac{\partial H}{\partial u} &= 2u+\lambda(1-\cos\theta)=0.\nonumber
\end{align} 
Thus, the optimal control $u$ satisfies
\begin{eqnarray}
	\label{eq:u*}
	u=-\frac{1}{2}\lambda (1-\cos\theta).
\end{eqnarray}
With \eqref{eq:u*} and \eqref{eq:lambda_dot}, this optimal control problem is transformed to a boundary value problem, whose solution  characterizes the optimal trajectories,
\begin{align}
	\label{eq:theta_p2}
    \dot{\theta} &=\alpha+\beta\cos\theta-\frac{\lambda}{2}(1-\cos\theta)^2,\\
	\label{eq:lambda_p2}
    \dot{\lambda} &= \lambda\beta+\frac{\lambda^2}{2}(1-\cos\theta)\sin\theta,
\end{align}
with boundary conditions $\theta(0)=0$ and $\theta(T)=2\pi$, while $\lambda_0=\lambda(0)$ and $\lambda(T)$ are unspecified.

Additionally, since the Hamiltonian is not explicitly dependent on time, the optimal triple $(\lambda,\theta,u)$ satisfies $H(\lambda,\theta,u)=c$, $\forall\, 0\leq t\leq T$,
where $c$ is a constant. Together with \eqref{eq:u*} and \eqref{eq:hamiltonian}, this yields
\begin{equation}
	\label{eq:quadratic}
	\lambda(\alpha+\beta\cos\theta)-\frac{\lambda^2}{4}(1-\cos\theta)^2=c.
\end{equation}
Since $\theta(0)=0$, $c=2\lambda_0$, where $\lambda_0$ is undetermined. The optimal multiplier can be found by solving the above quadratic equation \eqref{eq:quadratic}, which gives
\begin{equation}
	\label{eq:lambda}
\lambda=\frac{2(\alpha+\beta\cos\theta)\pm2\sqrt{(\alpha+\beta\cos\theta)^2-2\lambda_0(1-\cos\theta)}}{(1-\cos\theta)^2},
\end{equation}
and then, from \eqref{eq:theta_p2}, the optimal phase trajectory follows
\begin{equation}
	\label{eq:theta}
	\dot{\theta}=\mp\sqrt{(\alpha+\beta\cos\theta)^2-2\lambda_0(1-\cos\theta)}.
\end{equation}
Integrating \eqref{eq:theta}, we find the spiking time $T$ in terms of the initial condition $\lambda_0$,
\begin{equation}
	\label{eq:T2}
	T=\int_0^{2\pi}{\frac{1}{\sqrt{(\alpha+\beta\cos\theta)^2-2\lambda_0(1-\cos\theta)}}}d\theta.
\end{equation}
Note that we choose the positive sign in \eqref{eq:theta}, which corresponds to forward phase evolution. Therefore, given a desired spiking time $T$ of the neuron, the initial value $\lambda_0$ can be found via the one-to-one relation in \eqref{eq:T2}. Consequently, the optimal trajectories of $\theta$ and $\lambda$ can be easily computed by evolving \eqref{eq:theta_p2} and \eqref{eq:lambda_p2} forward in time. Plugging \eqref{eq:lambda} into \eqref{eq:u*},
we obtain the optimal feedback law for spiking a Theta neuron at time $T$,
\begin{equation}
\label{eq:u*2}
    u(t)^*= \frac{-(\alpha+\beta\cos\theta)+\sqrt{(\alpha+\beta\cos\theta)^2-2\lambda_0(1-\cos\theta)^2}}{1-\cos\t}.
\end{equation}
% where $\lambda_0$ is to be calculated according to \eqref{eq:T2}.

% ================== Theta Neuron (bounded control)==============
\subsubsection{Bounded Minimum-Power Control of a Theta Neuron}
\label{apd:thetaneuron_bounded}
Given the bound $M$ on the control amplitude, if $|u^*(t)|\leq M$ for all $t\in[0,T]$, then the amplitude constraint is inactive and obviously the bounded minimum-power control is given by \eqref{eq:u*2} and \eqref{eq:T2}. If, however, $|u^*(t)|> M$ for some time interval, e.g., $t\in[t_1,t_2]\subset[0,T]$, which alternatively corresponds to $|u^*(\t)|> M$ for $\t(t_1)=\t_1$, $\t(t_2)=\t_2$, and $\theta\in [\theta_1,\theta_2]\subset[0,2\pi]$, the amplitude constraint is active and the optimal control will depend on $M$. We first consider $u^*(\t)>M$ for $\t\in [\theta_1,\theta_2]$ and observe in this case that $u(\t)=M$ is the minimizer of the Hamiltonian $H$ as in \eqref{eq:hamiltonian}, since $H$ is convex in $u$. The Hamiltonian for this interval is then given by $H=M^2+\lambda(\alpha+\beta\cos\theta+M-M\cos\theta)$. Because, by the maximum principle, $H$ is a constant along the optimal trajectory, the Lagrange multiplier $\lambda$ is given by,
\begin{equation}
	\label{eq:lambda_bounded}
	\lambda=\frac{c-M^2}{\alpha +\beta\cos\theta +M-M\cos\theta},
\end{equation}
which satisfies the adjoint equation \eqref{eq:lambda_dot}. Therefore, $u(\t)=M$ is optimal for $\theta\in[\theta_1,\theta_2]$. The value of the constant $c=2\lambda_0$ can be determined by applying the initial conditions, $\t(0)=0$ and $\lambda(0)=\lambda_0$ to \eqref{eq:hamiltonian}. Similarly, we can show that $u(t)=-M$ is optimal when $u^*(\t)<-M$ for some $\t\in[\t_3,\t_4]\subset[0,2\pi]$. Consequently, the constrained optimal control can be synthesized according to \eqref{eq:bounded_control} and \eqref{eq:bounded_control_T}.

Note that the number of time intervals that $|u^*(\t)|>M$ defines the number of switches in the optimal control law. Specifically, if $|u^*(\t)|>M$ for $n$ time intervals, then the optimal control will have $2n$ switches. Here we consider the simplest case, where the optimal control has only two switches, which is actually the only case for the Theta neuron model. As a result, suppose that $u^*(\t)>M$ for only one time interval, and then there are two switching angles $\theta_1$ and $\theta_2$ at which $u^*(\theta_1)=u^*(\theta_2)=M$. These two conditions, together with \eqref{eq:bounded_control_T}, determine the unknown parameters $\theta_1$, $\theta_2$, and $\lambda_0$ that characterize the bounded optimal control, $u_M^*$, as given in \eqref{eq:bounded_control} for the specified spiking time $T$. Note that the range of feasible spiking times is determined by the bound of the control amplitude $M$. A complete characterization of possible spiking range can be found in \cite{Li_PRE11}.  

\bibliographystyle{IEEEtran}
\bibliography{prc_bib}

\begin{thebibliography}{10}
\providecommand{\url}[1]{#1}
\csname url@samestyle\endcsname
\providecommand{\newblock}{\relax}
\providecommand{\bibinfo}[2]{#2}
\providecommand{\BIBentrySTDinterwordspacing}{\spaceskip=0pt\relax}
\providecommand{\BIBentryALTinterwordstretchfactor}{4}
\providecommand{\BIBentryALTinterwordspacing}{\spaceskip=\fontdimen2\font plus
\BIBentryALTinterwordstretchfactor\fontdimen3\font minus
  \fontdimen4\font\relax}
\providecommand{\BIBforeignlanguage}[2]{{%
\expandafter\ifx\csname l@#1\endcsname\relax
\typeout{** WARNING: IEEEtran.bst: No hyphenation pattern has been}%
\typeout{** loaded for the language `#1'. Using the pattern for}%
\typeout{** the default language instead.}%
\else
\language=\csname l@#1\endcsname
\fi
#2}}
\providecommand{\BIBdecl}{\relax}
\BIBdecl

\bibitem{Strogatz01}
S.~Strogatz, \emph{Nonlinear Dynamics And Chaos: With Applications To Physics,
  Biology, Chemistry, And Engineering}, 1st~ed., ser. Studies in
  nonlinearity.\hskip 1em plus 0.5em minus 0.4em\relax Westview Press, 2001.

\bibitem{Uhlhaas06}
P.~Uhlhaas and W.~Singer, ``Neural synchrony in brain disorders: Relevance for
  cognitive dysfunctions and pathophysiology,'' \emph{Neuron}, vol.~52, no.~1,
  pp. 155--168, 2006.

\bibitem{Hanson78}
F.~Hanson, ``Comparative studies of firefly pacemakers,'' \emph{Federation
  proceedings}, vol.~38, no.~8, pp. 2158--2164, 1978.

\bibitem{Mirollo90}
R.~Mirollo and S.~Strogatz, ``Synchronization of pulse-coupled biological
  oscillators,'' \emph{SIAM Journal on Applied Mathematics}, vol.~50, no.~6,
  pp. 1645--1662, 1990.

\bibitem{Ermentrout84}
G.~Ermentrout and J.~Rinzel, ``Beyond a pacemaker's entrainment limit: phase
  walk-through,'' \emph{American Journal of Physiology - Regulatory,
  Integrative and Comparative Physiology}, vol. 246, no.~1, 1984.

\bibitem{Nishikawa08}
T.~Nishikawa, N.~Gulbahce, and A.~E. Motter, ``Spontaneous reaction silencing
  in metabolic optimization,'' \emph{PLoS Computational Biology}, vol.~4,
  no.~12, p. e1000236, 2008.

\bibitem{Fischer00}
I.~Fischer, Y.~Liu, and P.~Davis, ``Synchronization of chaotic semiconductor
  laser dynamics on subnanosecond time scales and its potential for chaos
  communication,'' \emph{Physical Review A}, vol.~62, 2000.

\bibitem{Blekhman88}
I.~Blekhman, \emph{Synchronization in science and technology}.\hskip 1em plus
  0.5em minus 0.4em\relax New York: ASME Press translations, 1988.

\bibitem{Harada10}
T.~Harada, H.~Tanaka, M.~Hankins, and I.~Kiss, ``Optimal waveform for the
  entrainment of a weakly forced oscillator,'' \emph{Physical Review Letters},
  vol. 105, no.~8, 2010.

\bibitem{Kiss07}
I.~Z. Kiss, C.~G. Rusin, H.~Kori, and J.~L. Hudson, ``Engineering complex
  dynamical structures: Sequential patterns and desynchronization,''
  \emph{Science}, vol. 316, no. 5833, pp. 1886--1889, 2007.

\bibitem{Ashwin92}
P.~Ashwin and J.~Swift, ``The dynamics of $n$ weakly coupled identical
  oscillators,'' \emph{Journal of Nonlinear Science}, vol.~2, no.~6, pp.
  69--108, 1992.

\bibitem{Benabid91}
A.~L. Benabid and P.~Pollak, ``Long-term suppression of tremor by chronic
  stimulation of the ventral intermediate thalamic nucleus,'' \emph{Lnacent},
  vol. 337, pp. 403--406, 1991.

\bibitem{Schiff94}
S.~Schiff, ``Controlling chaos in the brain,'' \emph{Nature}, vol. 370, pp.
  615--620, 1994.

\bibitem{Hoppensteadt01}
F.~C. Hoppensteadt and E.~M. Izhikevich, ``Synchronization of mems resonators
  and mechanical neurocomputing,'' \emph{IEEE Transactions On Circuits And
  Systems I-Fundamental Theory And Applications}, vol.~48, no.~2, pp. 133--138,
  2001.

\bibitem{Hoppensteadt00}
------, ``Synchronization of laser oscillators, associative memory, and optical
  neurocomputing,'' \emph{Physical Review E}, vol.~62, no.~3, pp. 4010--4013,
  2000.

\bibitem{Li_DSCC11}
A.~Zlotnik and J.-S. Li, ``Optimal asymptotic entrainment of phase-reduced
  oscillators,'' in \emph{ASME Dynamic Systems and Control Conference},
  Arlington, VA, October 2011.

\bibitem{Winfree80}
A.~T. Winfree, \emph{The Geometry of Biological Time}.\hskip 1em plus 0.5em
  minus 0.4em\relax Springer-Verlag, New York, 1980.

\bibitem{Glass91}
L.~Glass, ``Cardiac arrhythmias and circle maps-a classical problem,''
  \emph{Chaos}, vol.~1, no.~1, pp. 13--19, 1991.

\bibitem{Kuramoto84}
Y.~Kuramoto, \emph{Chemical Oscillations, Waves, and Turbulence}.\hskip 1em
  plus 0.5em minus 0.4em\relax New York: Springer, 1984.

\bibitem{Acebron05}
J.~A. Acebron, L.~L. Bonilla, C.~J.~P. Vicente, F.~Ritort, and R.~Spigler,
  ``The kuramoto model: A simple paradigm for synchronization phenomena,''
  \emph{Reviews Of Modern Physics}, vol.~77, no.~1, pp. 137--185, 2005.

\bibitem{Kiss05}
I.~Z. Kiss, Y.~M. Zhai, and J.~L. Hudson, ``Predicting mutual entrainment of
  oscillators with experiment-based phase models,'' \emph{Physical Review
  Letters}, vol.~94, no.~24, p. 248301, 2005.

\bibitem{Preyer05}
A.~J. Preyer and R.~J. Butera, ``Neuronal oscillators in aplysia californica
  that demonstrate weak coupling in vitro,'' \emph{Physical Review Letters},
  vol.~95, no.~13, p.~4, 2005.

\bibitem{Netoff05}
J.~C.~B. T.~I.~Netoff, C. D.~Acker and J.~A. White, ``Beyond two-cell networks:
  experimental measurement of neuronal responses to multiple synaptic inputs,''
  \emph{Journal Of Computational Neuroscience}, vol.~18, no.~3, pp. 287--295,
  2005.

\bibitem{Kano10}
T.~Kano and S.~Kinoshita, ``Control of individual phase relationship between
  coupled oscillators using multilinear feedback,'' \emph{Physical Review E},
  vol.~81, p. 026206, 2010.

\bibitem{Rusin10}
C.~G. Rusin, H.~Kori, I.~Z. Kiss, and J.~L. Hudson, ``Synchronization
  engineering: tuning the phase relationship between dissimilar oscillators
  using nonlinear feedback.'' \emph{Philosophical Transactions of the Royal
  Society A-Mathematical Physical And Engineering Sciences}, vol. 368, no.
  1918, pp. 2189--2204, 2010.

\bibitem{Zhai08}
Y.~Zhai, I.~Z. Kiss, and J.~L. Hudson, ``Control of complex dynamics with
  time-delayed feedback in populations of chemical oscillators:
  Desynchronization and clustering,'' \emph{Ind. \& Eng. Chem. Res.}, vol.~47,
  no.~10, pp. 3502--3514, 2008.

\bibitem{Li_PRE11}
I.~Dasanayake and J.-S. Li, ``Optimal design of minimum-power stimuli for phase
  models of neuron oscillators,'' \emph{Physical Review E}, vol.~83, p. 061916,
  2011.

\bibitem{Li_CDC11_Neuron}
------, ``Constrained minimum-power control of spiking neuron oscillators,'' in
  \emph{50th IEEE Conference on Decision and Control}, Orlando, FL, December
  2011.

\bibitem{Li_IEEE_CB}
------, ``Charge-balanced minimum-power controls for spiking neuron
  oscillators,'' \emph{IEEE Transactions on Automatic Control (under review)}.

\bibitem{Nabi09}
A.~Nabi and J.~Moehlis, ``Charge-balanced optimal input for phase models of
  spiking neurons,'' in \emph{Proc. ASME Dynamic System and Control
  Conference}, 2009, pp. 278--292.

\bibitem{Brown04}
E.~Brown, J.~Moehlis, and P.~Holmes, ``On the phase reduction and response
  dynamics of neural oscillator populations,'' \emph{Neural Computation},
  vol.~16, no.~4, pp. 673--715, 2004.

\bibitem{Izhikevich07}
E.~Izhikevich, \emph{Dynamical Systems in Neuroscience}, ser.
  Neuroscience.\hskip 1em plus 0.5em minus 0.4em\relax MIT Press, 2007.

\bibitem{Moehlis06}
J.~Moehlis, E.~Brown, and H.~Rabitz, ``Optimal inputs for phase models of
  spiking neurons,'' \emph{Journal of Computational and Nonlinear Dynamics},
  vol.~1, pp. 358--367, 2006.

\bibitem{ermentrout96}
B.~Ermentrout, ``Type \text{I} membranes, phase resetting curves, and
  synchrony,'' \emph{Neural Computation}, vol.~8, no.~5, pp. 979--1001, 1996.

\bibitem{Brockett76}
R.~Brockett, ``Nonlinear systems and differential geometry,'' \emph{Proceedings
  of the IEEE}, vol.~64, no.~1, pp. 61--72, 1976.

\bibitem{Levitan82}
B.~M. Levitan and V.~V. Zhikov, \emph{Almost periodic functions and
  differential equations}.\hskip 1em plus 0.5em minus 0.4em\relax Cambridge:
  Cambridge University Press, 1982.

\bibitem{Rose89}
R.~Rose and J.~Haindmarsh, ``The assembly of ionic currents in a thalamic
  neuron i. the three-dimensional model,'' \emph{Proc. R. Soc. Lond. B}, vol.
  237, pp. 267--28, 1989.

\bibitem{Keener98}
J.~Keener and J.~Sneyd, \emph{Mathematical Physiology}.\hskip 1em plus 0.5em
  minus 0.4em\relax New York: Springer-Verlag, 1998.

\bibitem{Sastry92}
S.~Sastry and R.~Montgomery, ``The structure of optimal controls for a steering
  problem,'' in \emph{Proc. IFAC Symposium on Nonlinear Control Systems}, 1992.

\bibitem{Montgomery89}
R.~Montgomery, ``Optimal control of deformable bodies, isoholonomic problems
  and sub-riemannian geometry,'' in \emph{Technical Report 05324-89,
  Mathematical Sciences Research Institute}, 1989.

\bibitem{Marks05}
W.~Marks, ``Deep brain stimulation for dystonia,'' \emph{Curr. Treat. Options
  Neurol.}, pp. 237--243, 2005.

\bibitem{Sussmann82}
H.~J. Sussmann, \emph{Time-Optimal Control in the Plane, in Feedback Control of
  Linear and Non-linear Systems, Lecture Notes in Control and Information
  Sciences}.\hskip 1em plus 0.5em minus 0.4em\relax Berlin: Springer-Verlag,
  1982.

\bibitem{Sussmann87}
------, ``The structure of time-optimal trajectories for single-input systems
  in the plane: The $c^{\infty}$ nonsingular case,'' \emph{SIAM J. Control
  Optim.}, vol.~25, pp. 433--465, 2004.

\bibitem{Li_SICON11}
D.~Stefanatos, H.~Schaettler, and J.-S. Li, ``Minimum-time frictionless atom
  cooling in harmonic traps,'' \emph{SIAM Journal on Control and Optimization},
  vol.~49, no.~6, pp. 2440--2462, 2011.

\bibitem{Isidori95}
A.~Isidodri, \emph{Nonlinear Control Systems}.\hskip 1em plus 0.5em minus
  0.4em\relax London: Springer, 1995.

\bibitem{Li_PNAS11}
J.-S. Li, J.~Ruths, T.-Y. Yu, H.~Arthanari, and G.~Wagner, ``Optimal pulse
  design in quantum control: A unified computational method,''
  \emph{Proceedings of the National Academy of Sciences}, vol. 108, no.~5, pp.
  1879--1884, 2011.

\bibitem{Ruths_JCP11}
J.~Ruths and J.-S. Li, ``A multidimensional pseudospectral method for optimal
  control of quantum ensembles,'' \emph{Journal of Chemical Physics}, vol. 134,
  p. 044128, 2011.

\bibitem{Li_IEEE_QCP11}
------, ``Optimal control of inhomogeneous ensembles,'' \emph{IEEE Transactions
  on Automatic Control: Special Issue on Control of Quantum Mechanical Systems
  (in press)}.

\bibitem{Li_CDC11_PS}
J.~Ruths, A.~Zlotnik, and J.-S. Li, ``Convergence of a pseudospectral method
  for optimal control of complex dynamical systems,'' in \emph{50th IEEE
  Conference on Decision and Control}, Orlando, FL, December 2011.

\bibitem{Canuto06}
C.~Canuto, M.~Y. Hussaini, A.~Quarteroni, and T.~A. Zang, \emph{Spectral
  Methods}.\hskip 1em plus 0.5em minus 0.4em\relax Berlin: Springer, 2006.

\bibitem{Elnagar95}
G.~Elnagar, M.~A. Kazemi, and M.~Razzaghi, ``The pseudospectral legendre method
  for discretizing optimal control problems,'' \emph{IEEE Transactions on
  Automatic Control}, vol.~40, no.~10, pp. 1793--1796, 1995.

\bibitem{Ross03}
M.~Ross and F.~Fahroo, ``Legendre pseudospectral approximations of optimal
  control problems,'' in \emph{New Trends in Nonlinear Dynamics and Control},
  W.~Kang, M.~Xiao, and C.~R. Borges, Eds.\hskip 1em plus 0.5em minus
  0.4em\relax Berlin: Springer, 2003, pp. 327--342.

\bibitem{Fahroo01}
F.~Fahroo and I.~Ross, ``Costate estimation by a legendre pseudospectral
  method,'' \emph{Journal of Guidance, Control, and Dynamics}, vol.~24, no.~2,
  pp. 270--277, 2001.

\bibitem{Fornberg98}
B.~Fornberg, \emph{A Practical Guide to Pseudospectral Methods}.\hskip 1em plus
  0.5em minus 0.4em\relax Cambridge University Press, 1998.

\bibitem{Szego59}
G.~Szego, \emph{Orthogonal Polynomials}.\hskip 1em plus 0.5em minus 0.4em\relax
  New York: American Mathematical Society, 1959.

\end{thebibliography}
\end{document}